# A NOTE ON THE APPLICATION OF THE GUERMOND-PASQUETTI MASS LUMPING CORRECTION TECHNIQUE FOR CONVECTION-DIFFUSION PROBLEMS

SERGII V. SIRYK

National Technical University of Ukraine "Igor Sikorsky Kyiv Polytechnic Institute", Kiev, Ukraine
*(Prospect Peremogy 37, Academic building 15, Kyiv 03056, Ukraine)*

**Abstract:** We provide a careful Fourier analysis of the Guermond-Pasquetti mass lumping correction technique *[Guermond J.-L., Pasquetti R. A correction technique for the dispersive effects of mass lumping for transport problems // Computer Methods in Applied Mechanics and Engineering. – 2013. – Vol. 253. – P. 186-198]* applied to pure transport and convection-diffusion problems. In particular, it is found that increasing the number of corrections reduces the accuracy for problems with diffusion; however all the corrected schemes are more accurate than the consistent Galerkin formulation in this case. For the pure transport problems the situation is the opposite. We also investigate the differences between two numerical solutions – the consistent solution and the corrected ones, and show that increasing the number of corrections makes solutions of the corrected schemes closer to the consistent solution in all cases.

## 1. Introduction

Applying the finite element methods for the spatial approximation of a non-stationary convection-diffusion equation usually yields a semi-discrete problem [1–7], i.e., a system of ordinary differential equations (SODE) of the form $M\dot{\vec{a}} = \vec{F}(t,\vec{a})$, where $\vec{a}(t)$ is the vector of unknown coefficients of the expansion of the numerical solution with respect to the trial functions, $\vec{F}$ is some vector function, $M$ is the so-called (consistent) mass matrix [1–4], which is sparse, non-diagonal and even non-symmetric in the general case. In the subsequent passage from SODE to difference schemes (when time derivatives are replaced by differences) the resulting difference schemes become implicit due to the non-diagonality of the matrix $M$. In addition, the matrix $M$ turns out to be time-dependent in some statements of numerical problems, which can lead to the necessity to inverse/factorize the matrix $M$ at each time step of integration of the SODE obtained [4–7]. The mass lumping technique [1–4, 7–11] is often used in computational practice to facilitate computational efforts and avoid the necessity for sophisticated (and computationally high-cost) algebra. The essence of this technique is to replace the matrix $M$ by a diagonal matrix (we denote it by $\overline{M}$). There are several options for implementing this technique [3, 9–11]. Since in what follows we analyze only the case with continuous piecewise-linear (Lagrange) trial functions, we simply use the sums of elements in the corresponding rows of the matrix $M$ to obtain the diagonal elements of the matrix $\overline{M}$, which is the standard and commonly used procedure in the literature (see [1–3, 7–11]), also equivalent to choosing the interpolation points as quadrature points for numerical integration in this case [11]. It should be noted that this operation may produce non-positive definite (with zero or negative diagonal elements) mass matrices for higher-order elements [3, 10, 11], therefore special attention (and special constructions like quasi-lumping introduced in [11]) may be needed in this case. After performing mass lumping we obtain the "lumped" SODE $\overline{M}\dot{\vec{a}} = \vec{F}(t,\vec{a})$ instead of the original "consistent" SODE. The use of mass lumping makes it possible to treat the partial derivative with respect to time in the finite-element method (FEM) schemes in the same way as it is done in the finite difference methods. It is clear that there is no need to perform time-consuming inversion operations after carrying out this diagonalization of the matrix $M$. Note also that mass lumping plays an important role in the construction of maximum principle preserving methods (see [12–14] and the references therein).

However, it is known (see [3, 7–11]) that the application of mass lumping can introduce dispersion and dissipation errors into numerical schemes leading to significant inaccuracies in the numerical solution (a detailed review of the "pro" and "contra" papers of using mass lumping in numerical schemes is given in [8]; see also [4] for a review of avoiding the inversion of mass matrix). The fundamental paper of Guermond J.-L. and Pasquetti R. [11] is devoted to the investigation and overcoming this essential drawback of mass lumping. Their technique is based on using a matrix series to approximate the matrix $M^{-1}$: the matrix $M$ is represented in the form $M = \overline{M}(I - A)$ where $I$ is the identity matrix, the matrix $A \equiv \overline{M}^{-1}(\overline{M} - M)$, from which we get $M^{-1} = (I + A + A^2 + A^3 + \ldots)\overline{M}^{-1}$ (Neumann series). It should be noted that only the pure (diffusionless) transport equation was considered in [11], and the main focus there was on semi-discrete approximations of the standard classical Galerkin method (without stabilization) – in particular, the convergence of the corresponding Neumann matrix series was rigorously proved for the classical Galerkin FEM with linear elements. The authors note that the use of even one (the first) correction term (i.e., $(I + A)\overline{M}^{-1}$ instead of $M^{-1}$) can significantly improve the accuracy of the numerical solution, and that correcting the lumped mass matrix four times ($M^{-1} \approx (I + A + A^2 + A^3 + A^4)\overline{M}^{-1}$) makes numerical results practically indistinguishable from those obtained by the consistent formulation. However, the authors did not make any detailed theoretical estimates of the quality and accuracy of the solution depending on the number of terms of the matrix series taken. Despite the computational attractiveness and wide application of this powerful technique (e.g., in constructing maximum principle preserving methods [12–14], level set methods for two-phase flows [15], various engineering applications [16, 17] etc.), and the overall high quotability of the paper [11], these important issues remain still unexplored. Therefore, our paper can be considered as the first attempt to clarify these issues.

Our paper provides a careful Fourier analysis of this technique in application to pure transport and convection-diffusion problems. We show that increasing the number of corrections (i.e., the number of correcting terms in the Neumann series) leads to error increase in the presence of diffusion terms (see details in Lemma 1 and Proposition 3 below) – thus, in contrast to pure transport problems, it is advisable to use only one (the first) correction for problems with diffusion. We also show that all the corrected schemes are more accurate than the consistent Galerkin formulation for problems with diffusion (see details in Lemma 2 and Proposition 5 below). These results seem to be new, quite unexpected, and unnoticed earlier in the literature. For the pure (diffusionless) transport problems the situation is completely opposite – i.e., increasing the number of corrections improves the accuracy of the numerical solution (see details in Proposition 4 below), and the consistent Galerkin formulation produces more accurate results than all the corrected schemes (see details in Proposition 6 below). We also investigate the differences between two numerical solutions – the consistent solution and the corrected ones, and show that increasing the number of corrections makes solutions of the corrected schemes closer to the consistent solution in all cases (see details in Proposition 8 below). Despite the one-dimensionality of the presented Fourier analysis (but it can be extended to tensor product meshes in higher dimensions), the numerical examples given below show that, in general, one may conjecture the validity of such results in any dimensions including unstructured non-uniform meshes (with randomly distributed nodes).

## 2. Fourier analysis of semi-discrete approximations

### 2.1. Continuous and numerical problems setting

We consider the one-dimensional convection-diffusion equation [1–3, 5, 7, 18–20]

$$\frac{\partial u}{\partial t} + \lambda \frac{\partial u}{\partial x} - \kappa \frac{\partial^2 u}{\partial x^2} = 0, \tag{1}$$

where the real coefficients $\lambda$ and $\kappa$ are given, $\kappa \geq 0$, $\kappa + |\lambda| > 0$, and $u = u(t, x)$ is the unknown solution. To facilitate the use of Fourier analysis [1, 2, 18–19] we assume that the spatial mesh is uniform with mesh size $h > 0$ and mesh nodes $x_k = kh$, and that the coefficients $\lambda$ and $\kappa$ are constant. Using the standard continuous piecewise-linear trial functions (Lagrange functions, or "hat" functions) in the classical Galerkin finite-element formulation for (1), we obtain the SODE where the typical $k$-th equation (corresponding to the typical mesh node $x_k$) has the form (see [1, 2, 7, 18–19])

$$\frac{h}{6}\frac{da_{k-1}}{dt} + \frac{2h}{3}\frac{da_k}{dt} + \frac{h}{6}\frac{da_{k+1}}{dt} + \lambda \frac{a_{k+1} - a_{k-1}}{2} - \kappa \frac{a_{k+1} - 2a_k + a_{k-1}}{h} = 0. \tag{2}$$

Here $\vec{a} = \{a_k(t)\}$ are the coefficients of the expansion of the approximate solution with respect to the corresponding trial functions. Applying mass lumping to (2), we obtain the following equation:

$$h\frac{da_k}{dt} + \lambda \frac{a_{k+1} - a_{k-1}}{2} - \kappa \frac{a_{k+1} - 2a_k + a_{k-1}}{h} = 0. \tag{3}$$

Thus, mass matrices $M$ and $\overline{M}$ have the following obvious representations:

$$M_{ij} = \begin{cases} h/6, & \text{if } j = i \pm 1, \\ 2h/3, & \text{if } j = i, \\ 0, & \text{otherwise,} \end{cases} \qquad \overline{M}_{ij} = \begin{cases} h, & \text{if } j = i, \\ 0, & \text{otherwise} \end{cases}$$

for the typical mesh node $x_i$.

Equations (2) and (3) form the consistent and lumped semi-discrete Galerkin FEM formulations (with linear elements), respectively.

Introducing finite-difference operators $C$ and $D$ of the central first and second derivatives with the matrix representations

$$C_{ij} = \begin{cases} \pm 1/(2h), & \text{if } j = i \pm 1, \\ 0, & \text{otherwise,} \end{cases} \qquad D_{ij} = \begin{cases} 1/h^2, & \text{if } j = i \pm 1, \\ -2/h^2, & \text{if } j = i, \\ 0, & \text{otherwise} \end{cases}$$

for the typical mesh node $x_i$, systems (2) and (3) can be rewritten as $M\dot{\vec{a}} + \lambda h \cdot C\vec{a} - \kappa h \cdot D\vec{a} = 0$ and $\overline{M}\dot{\vec{a}} + \lambda h \cdot C\vec{a} - \kappa h \cdot D\vec{a} = 0$, respectively.

For the standard issues of setting and handling the initial/boundary conditions in Galerkin FEM formulations one can see for example [1–3, 11, 18–20].

### 2.2. Mass lumping corrections

The generic form of the consistent SODE $M\dot{\vec{a}} = \vec{F}(t, \vec{a})$ can be rewritten as $\dot{\vec{a}} = M^{-1}\vec{F}(t, \vec{a})$. Approximating the matrix $M^{-1}$ in the manner described above we arrive to the following definition.

**Definition 1.** The system $\dot{\vec{a}} = (I + A + A^2 + A^3 + \ldots + A^n)\overline{M}^{-1}\vec{F}(t, \vec{a})$ is called the $n$-th corrected scheme.

Note that from this definition we obtain the standard (non-corrected) lumped semi-discrete scheme $\vec{M}\dot{\vec{a}} = \vec{F}(t, \vec{a})$ for $n = 0$.

Let us introduce the following notation (see [21]): $y_{\bar{x},k} = (y_k - y_{k-1})/h$, $y_{x,k} = (y_{k+1} - y_k)/h$,

$y_k^{(1)} = y_{\tilde{x},k} = (y_{\bar{x},k} + y_{x,k})/2 = (y_{k+1} - y_{k-1})/(2h)$, $\quad y_k^{(2)} = y_{\bar{x}x,k} = (y_{\bar{x}})_{x,k} = (y_{k+1} - 2y_k + y_{k-1})/h^2$,

$y_k^{(3)} = y_{\tilde{x}\bar{x}x,k} = (y_{\tilde{x}})_{\bar{x}x,k} = (y_{k+2} - 2y_{k+1} + 2y_{k-1} - y_{k-2})/(2h^3)$,

$y_k^{(4)} = y_{\bar{x}x\bar{x}x,k} = (y_{\bar{x}x})_{\bar{x}x,k} = (y_{k+2} - 4y_{k+1} + 6y_k - 4y_{k-1} + y_{k-2})/h^4$,

$y_k^{(5)} = y_{\tilde{x}\bar{x}x\bar{x}x,k} = (y_{\tilde{x}})_{\bar{x}x\bar{x}x,k} = (y_{k+3} - 4y_{k+2} + 5y_{k+1} - 5y_{k-1} + 4y_{k-2} - y_{k-3})/(2h^5)$,

$y_k^{(6)} = y_{\bar{x}x\bar{x}x\bar{x}x,k} = (y_{\bar{x}x})_{\bar{x}x\bar{x}x,k} = (y_{k+3} - 6y_{k+2} + 15y_{k+1} - 20y_k + 15y_{k-1} - 6y_{k-2} + y_{k-3})/h^6$, ... ,

etc. With usage of the operators $C$ and $D$ we can also represent these derivatives as $y_k^{(2n)} = (D^n \vec{y})_k$, $y_k^{(2n+1)} = (D^n C \vec{y})_k$, $n \geq 0$, for arbitrary mesh function $\vec{y}$ (see [21]).

Note that for the case under consideration we have $\vec{F} = \kappa h \cdot D\vec{a} - \lambda h \cdot C\vec{a}$, and corresponding corrected schemes are characterized by the following.

**Proposition 1.** The typical $k$-th equation (corresponding to the typical mesh node $x_k$) of the $n$-th corrected scheme has the form

$$\frac{da_k}{dt} + \lambda(C\vec{a})_k - \kappa(D\vec{a})_k + \sum_{m=1}^{n}(-1)^m \frac{h^{2m}}{6^m}\left(\lambda(D^m C\vec{a})_k - \kappa(D^{m+1}\vec{a})_k\right) = 0, \quad (4.1)$$

or, using the introduced notation, it can be rewritten as

$$\frac{da_k}{dt} + \lambda a_k^{(1)} - \kappa a_k^{(2)} - \frac{\lambda h^2}{6}a_k^{(3)} + \frac{\kappa h^2}{6}a_k^{(4)} + \frac{\lambda h^4}{36}a_k^{(5)} - \frac{\kappa h^4}{36}a_k^{(6)} + \\ + \ldots + (-1)^n \frac{\lambda h^{2n}}{6^n}a_k^{(2n+1)} - (-1)^n \frac{\kappa h^{2n}}{6^n}a_k^{(2n+2)} = 0. \quad (4.2)$$

**Proof.** A direct calculation shows that the matrix $A$ has the following matrix representation:

$$A_{ij} = \begin{cases} -1/6, & \text{if } j = i \pm 1, \\ 1/3, & \text{if } j = i, \\ 0, & \text{otherwise,} \end{cases}$$

thereby $(A\vec{y})_k = (-h^2/6)y_{\bar{x}x,k} = (-h^2/6)(D\vec{y})_k$ and $(A^n \vec{y})_k = ((-1)^n h^{2n}/6^n)(D^n \vec{y})_k$ for arbitrary mesh function $\vec{y}$. The last equation and the representation $\vec{F} = \kappa h \cdot D\vec{a} - \lambda h \cdot C\vec{a}$ imply (4.1)-(4.2). ∎

*2.3. Fourier analysis*

Due to the Fourier approach [1, 2, 18–20, 22, 23] we are looking for particular solutions of the semi-discrete approximation (2) in the form of harmonics $a_k(t) = e^{\omega t} e^{ikhp}$, where $\omega$ is some complex number to be determined, $i^2 = -1$, $p$ is the real number (the spatial wave number of the harmonic). We also assume that $|ph| \leq \pi$, which is standard condition in Fourier analysis (related to the Nyquist limit) [1, 18–20, 22, 23]. Substituting $a_k(t) = e^{\omega t} e^{ikhp}$ into (2), we obtain the following expression for $\omega$ (we denote it by $\omega_G$ for distinguishing):

$$\omega_G = \frac{-\dfrac{4\kappa}{h^2}\sin^2\left(\dfrac{ph}{2}\right) - i\dfrac{\lambda}{h}\sin(ph)}{1 - \dfrac{2}{3}\sin^2\left(\dfrac{ph}{2}\right)}.$$

Let us denote the number $\omega$ for the $n$-th corrected scheme (equations (4.1)-(4.2)) by $\omega_n$.

**Proposition 2.** The real and imaginary parts of $\omega_n$ can be expressed as follows:

$$\mathrm{Re}(\omega_n) = -\frac{4\kappa}{h^2}\sin^2\left(\frac{ph}{2}\right) - \frac{8\kappa}{3h^2}\sin^4\left(\frac{ph}{2}\right) - \frac{16\kappa}{9h^2}\sin^6\left(\frac{ph}{2}\right) - \ldots - \frac{2^{n+2}\kappa}{3^n h^2}\sin^{2n+2}\left(\frac{ph}{2}\right),$$

$$\mathrm{Im}(\omega_n) = -\frac{\lambda}{h}\sin(ph) - \frac{2\lambda}{3h}\sin(ph)\sin^2\left(\frac{ph}{2}\right) - \frac{4\lambda}{9h}\sin(ph)\sin^4\left(\frac{ph}{2}\right) - \ldots - \frac{2^n\lambda}{3^n h}\sin(ph)\sin^{2n}\left(\frac{ph}{2}\right).$$

**Proof.** Substituting $a_k(t) = e^{\omega_n t}e^{ikhp}$ into the equation (4.2), using Euler formulas (and, in particular, the relation $e^{iph} - 2 + e^{-iph} = -4\sin^2(ph/2)$) and the representation of the coefficients of central differences via binomial coefficients (e.g., see [24]), after arithmetic transformations we obtain

$$\omega_n = \frac{\kappa}{h^2}\left(-4\sin^2\left(\frac{ph}{2}\right)\right) - \frac{\lambda}{2h}2i\sin(ph) - \frac{\kappa h^2}{6}\frac{1}{h^4}\left(-4\sin^2\left(\frac{ph}{2}\right)\right)^2 + \frac{\lambda h^2}{6}\frac{1}{2h^3}2i\sin(ph)\left(-4\sin^2\left(\frac{ph}{2}\right)\right) +$$

$$+ \ldots + (-1)^n \frac{\kappa h^{2n}}{6^n}\frac{1}{h^{2n+2}}\left(-4\sin^2\left(\frac{ph}{2}\right)\right)^{n+1} - (-1)^n \frac{\lambda h^{2n}}{6^n}\frac{1}{2h^{2n+1}}2i\sin(ph)\left(-4\sin^2\left(\frac{ph}{2}\right)\right)^n,$$

which is equivalent to the above expressions. ∎

Substituting the ansatz $u(t,x) = e^{\overline{\omega}t}e^{ixp}$ into the equation (1) we obtain $\overline{\omega} = -\kappa p^2 - i\lambda p$. Note that $\lim_{h\to 0}\omega_G = \lim_{h\to 0}\omega_n = \overline{\omega}$, i.e., the solutions (sought in the class of harmonics) of all considered numerical problems tend to the solution of the problem (1) if $h$ tends to zero.

**Lemma 1.** Let $\kappa > 0$, $n \geq 1$, and let $\lambda$ be arbitrary. Then for any $p \neq 0$ the representation $|\omega_{n+1} - \overline{\omega}|^2 - |\omega_n - \overline{\omega}|^2 = \left(\dfrac{2^{n+1}\kappa^2 p^4}{3^{n+1}}\sin^{2n+2}\left(\dfrac{z}{2}\right)\right)\cdot f_n(z)$ is true, where $z \equiv ph$ and $f_n(z)$ is defined by

$$f_n(z) = \frac{4\sin^2(z/2)}{z^2}\left(\frac{8\sin^2\dfrac{z}{2}}{z^2}\cdot\frac{1 - \left(\dfrac{2}{3}\sin^2\dfrac{z}{2}\right)^{n+2}}{1 - \dfrac{2}{3}\sin^2\dfrac{z}{2}} - 2 - \frac{2^{n+3}}{3^{n+1}z^2}\sin^{2n+4}\frac{z}{2}\right) +$$

$$+ \left(\frac{\lambda}{\kappa}\right)^2 \frac{\sin z}{p^2 z}\left(2\frac{\sin z}{z}\cdot\frac{1 - \left(\dfrac{2}{3}\sin^2\dfrac{z}{2}\right)^{n+2}}{1 - \dfrac{2}{3}\sin^2\dfrac{z}{2}} - 2 - \frac{2^{n+1}\sin z}{3^{n+1}z}\sin^{2n+2}\frac{z}{2}\right).$$

Moreover, $f_n(z) \geq f_{n-1}(z) \geq \ldots \geq f_1(z)$ for arbitrary $z \neq 0$ (these inequalities are strict if $z$ is not a multiple of $2\pi$), and $f_1(z) > 0$ for arbitrary $z$ satisfying $0 < |z| \leq z_0$, where

$$z_0 = \sqrt{\frac{6\left(98(\lambda/(\kappa p))^2 + 91 - \sqrt{10661 - 4004(\lambda/(\kappa p))^2 + 9604(\lambda/(\kappa p))^4}\right)}{156(\lambda/(\kappa p))^2 - 17}}$$ in the case $156(\lambda/\kappa)^2 \neq 17p^2$

and $z_0 = 6\sqrt{130/1133}$ in the case $156(\lambda/\kappa)^2 = 17p^2$. The real number $z_0$ defined above always exists, positive and is less than $\pi$.

**Proof.** Let us denote $A_n = \dfrac{2^{n+2}\kappa}{3^n h^2}\sin^{2n+2}\left(\dfrac{ph}{2}\right)$, $B_n = \dfrac{2^n \lambda}{3^n h}\sin(ph)\sin^{2n}\left(\dfrac{ph}{2}\right)$. Then we have

$$|\omega_{n+1} - \overline{\omega}|^2 - |\omega_n - \overline{\omega}|^2 = |(\omega_n - A_{n+1} - iB_{n+1}) - \overline{\omega}|^2 - |\omega_n - \overline{\omega}|^2 = (\text{Re}(\omega_n - \overline{\omega}) - A_{n+1})^2 +$$

$$+ (\text{Im}(\omega_n - \overline{\omega}) - B_{n+1})^2 - (\text{Re}(\omega_n - \overline{\omega}))^2 - (\text{Im}(\omega_n - \overline{\omega}))^2 = A_{n+1}(A_{n+1} - 2\text{Re}(\omega_n - \overline{\omega})) +$$

$$+ B_{n+1}(B_{n+1} - 2\text{Im}(\omega_n - \overline{\omega})) = B_{n+1}\left(B_{n+1} + 2\dfrac{\lambda\sin(ph)}{h}\left(1 + \dfrac{2}{3}\sin^2\dfrac{ph}{2} + \ldots + \dfrac{2^n}{3^n}\sin^{2n}\dfrac{ph}{2}\right) - 2\lambda p\right) +$$

$$+ A_{n+1}\left(A_{n+1} + 2\left(\dfrac{4\kappa}{h^2}\sin^2\dfrac{ph}{2} + \dfrac{8\kappa}{3h^2}\sin^4\dfrac{ph}{2} + \ldots + \dfrac{2^{n+2}\kappa}{3^n h^2}\sin^{2n+2}\dfrac{ph}{2} - \kappa p^2\right)\right) =$$

$$= B_{n+1}\left(2\dfrac{\lambda\sin(ph)}{h}\left(1 + \dfrac{2}{3}\sin^2\dfrac{ph}{2} + \ldots + \dfrac{2^{n+1}}{3^{n+1}}\sin^{2(n+1)}\dfrac{ph}{2}\right) - 2\lambda p - \dfrac{2^{n+1}\lambda}{3^{n+1}h}\sin(ph)\sin^{2n+2}\dfrac{ph}{2}\right) +$$

$$+ A_{n+1}\left(\dfrac{8\kappa}{h^2}\sin^2\dfrac{ph}{2}\left(1 + \dfrac{2}{3}\sin^2\dfrac{ph}{2} + \ldots + \dfrac{2^{n+1}}{3^{n+1}}\sin^{2(n+1)}\dfrac{ph}{2}\right) - 2\kappa p^2 - \dfrac{2^{n+3}\kappa}{3^{n+1}h^2}\sin^{2n+4}\dfrac{ph}{2}\right);$$

collapsing the expressions in the inner parentheses (the geometric progression with the common ratio $(2/3)\sin^2(ph/2)$), after arithmetic transformations we obtain the representation

$$|\omega_{n+1} - \overline{\omega}|^2 - |\omega_n - \overline{\omega}|^2 = \left(\dfrac{2^{n+1}\kappa^2 p^4}{3^{n+1}}\sin^{2n+2}\left(\dfrac{z}{2}\right)\right)\cdot f_n(z).$$

Next, $f_{n+1}(z) - f_n(z) = \left(1 + \dfrac{2}{3}\sin^2\dfrac{z}{2}\right)\dfrac{2^{n+5}}{3^{n+1}z^4}\sin^{2n+6}\dfrac{z}{2} + \dfrac{(\lambda/\kappa)^2\sin^2 z}{p^2 z^2}\left(1 + \dfrac{2}{3}\sin^2\dfrac{z}{2}\right)\dfrac{2^{n+1}}{3^{n+1}}\sin^{2n+2}\dfrac{z}{2} \geq 0$.

Finally, let us show that $f_1(z) > 0$ for all $z$ satisfying $0 < |z| \leq z_0$, where $z_0$ is defined above. Using standard double-angle cosine formulas, we first rewrite the expression for $f_1(z)$ in the form

$$f_1(z) = \dfrac{371 - 524\cos z + 172\cos 2z - 20\cos 3z + \cos 4z - 72z^2(1 - \cos z)}{18z^4} +$$

$$+ \dfrac{(\lambda/\kappa)^2(101 - 16\cos z - 100\cos 2z + 16\cos 3z - \cos 4z - 144z\sin z)}{72p^2 z^2},$$

and then expand trigonometric expressions into the Taylor series with respect to $z$, which provides

$$f_1(z) = \dfrac{z^2}{30240}\left(5040 - 12z^2(98\mu^2 + 91) - (17 - 156\mu^2)z^4\right) + \sum_{m=4}^{+\infty}(-1)^m c_m z^{2m}, \quad \text{where} \quad \mu \equiv \lambda/(\kappa p),$$

$$c_m = \dfrac{172\cdot 2^{2m+4} - 20\cdot 3^{2m+4} + 4^{2m+4} - 524}{18\cdot(2m+4)!} + \mu^2\dfrac{4(1 - 3^{2m+2}) + 25\cdot 2^{2m+2} + 4^{2m+1}}{18\cdot(2m+2)!} - \dfrac{4(1 + (m+1)\mu^2)}{(2m+2)!}.$$

A direct arithmetic calculation using mathematical induction shows that $c_m > 0$ for arbitrary $m \geq 4$ and $\mu$, and that the sequence $\{c_m/c_{m+1}\}_{m\geq 4}$ increases monotonically for increasing $m$, which ensures that

$c_m z^{2m} > c_{m+1}z^{2(m+1)}$ ($\Leftrightarrow z^2 < c_m/c_{m+1}$) for all $m \geq 4$ if $0 < |z| < z_* \equiv \sqrt{\dfrac{c_4}{c_5}} = \dfrac{6\sqrt{33(80\kappa^2 p^2 + 79\lambda^2)}}{\sqrt{13399\kappa^2 p^2 + 30146\lambda^2}}$.

Note also that $c_m \xrightarrow[m \to \infty]{} 0$ and $z^{2m} c_m \xrightarrow[m \to \infty]{} 0$ for arbitrary fixed $z$. Thus, from the alternating series test (Leibniz rule) [25] we obtain that $0 < \sum_{m=4}^{+\infty}(-1)^m c_m z^{2m} \leq c_4 z^8$ and thereby $f_1(z) > z^2(5040 - 12z^2(98\mu^2 + 91) + z^4(156\mu^2 - 17))/30240$ when $0 < |z| < z_*$. The right part of the last expression is strictly positive for $0 < |z| < z_0$, where $z_0$ is the smallest positive root of the equation $5040 - 12z^2(98\mu^2 + 91) + z^4(156\mu^2 - 17) = 0$. It is easy to see that $z_0$ is determined by the expression given in the statement of this lemma. A direct arithmetic calculation shows that $z_0 < z_* < \pi$. Note that $z_0$ always exists, since the expression under the external square root and the discriminant of the last equation (which is equal to $144(98\mu^2 + 91)^2 - 4 \cdot 5040 \cdot (156\mu^2 - 17) = 144((686\mu^2 - 143)^2 + 501940)/49$) are both positive. ∎

**Remark 1.** Let $\kappa = 0$, $n \geq 1$, and $\lambda \neq 0$. Then we get the following equality (where $z \equiv ph \neq 0$):

$$|\omega_{n+1} - \overline{\omega}|^2 - |\omega_n - \overline{\omega}|^2 = \left(\frac{2^{n+1}\lambda^2 p^2}{3^{n+1}} \sin^{2n+2}\frac{z}{2}\right) \tilde{f}_n(z),$$ 

where the function $\tilde{f}_n(z)$ is defined by

$$\tilde{f}_n(z) = \frac{\sin z}{z}\left(\frac{\sin z}{z}\left(2 \cdot \frac{1 - \left(\frac{2}{3}\sin^2\frac{z}{2}\right)^{n+2}}{1 - \frac{2}{3}\sin^2\frac{z}{2}} - \frac{2^{n+1}}{3^{n+1}}\sin^{2n+2}\frac{z}{2}\right) - 2\right).$$ 

Moreover, $\tilde{f}_n(z) \geq \tilde{f}_{n-1}(z) \geq \ldots \geq \tilde{f}_1(z)$ and $\tilde{f}_n(z) \leq 2\frac{\sin z}{z}\left(\frac{\sin z}{z}\left(1 - \frac{2}{3}\sin^2\frac{z}{2}\right)^{-1} - 1\right)$ for arbitrary $n \geq 1$ and $z \neq 0$.

**Proof.** Following the beginning of the proof of Lemma 1, it is easy to see that in this case we have $|\omega_{n+1} - \overline{\omega}|^2 - |\omega_n - \overline{\omega}|^2 = B_{n+1}(B_{n+1} - 2\,\text{Im}(\omega_n - \overline{\omega}))$. The last expression is equal to the expression

$$\left(2\frac{\sin^2 z}{z^2}\left(1 + \frac{2}{3}\sin^2\frac{z}{2} + \ldots + \frac{2^{n+1}}{3^{n+1}}\sin^{2n+2}\frac{z}{2}\right) - 2\frac{\sin z}{z} - \frac{2^{n+1}}{3^{n+1}}\frac{\sin^2 z}{z^2}\sin^{2n+2}\frac{z}{2}\right)\frac{2^{n+1}\lambda^2 p^2}{3^{n+1}}\sin^{2n+2}\frac{z}{2},$$

which is equivalent to the above representation. Next,

$$\tilde{f}_{n+1}(z) - \tilde{f}_n(z) = \frac{\sin^2 z}{z^2}\left(1 + \frac{2}{3}\sin^2\frac{z}{2}\right)\frac{2^{n+1}}{3^{n+1}}\sin^{2n+2}\frac{z}{2} \geq 0, \quad 2\frac{\sin z}{z}\left(\frac{\sin z}{z}\left(1 - \frac{2}{3}\sin^2\frac{z}{2}\right)^{-1} - 1\right) - \tilde{f}_n(z) =$$

$$= \frac{\sin^2 z}{z^2}\left(\frac{2^{n+1}}{3^{n+1}}\sin^{2n+2}\frac{z}{2} + 2\sum_{m=n+2}^{+\infty}\frac{2^m}{3^m}\sin^{2m}\frac{z}{2}\right) \geq 0. \blacksquare$$

**Lemma 2.** Let $\kappa > 0$, $n \geq 1$, and let $\lambda$ be arbitrary. Then for any $p \neq 0$ the representation

$$|\omega_G - \overline{\omega}|^2 - |\omega_n - \overline{\omega}|^2 = \left(\frac{2^{n+1}\kappa^2 p^4}{3^{n+1}}\left(1 - \frac{2}{3}\sin^2\frac{z}{2}\right)^{-1}\sin^{2n+2}\frac{z}{2}\right) g_n(z)$$ 

is true, where $z \equiv ph$ and

$$g_n(z) = \frac{4}{z^2}\sin^2\left(\frac{z}{2}\right)\left(\frac{8}{z^2}\sin^2\frac{z}{2}\left(1-\frac{2}{3}\sin^2\frac{z}{2}\right)^{-1} - 2 - \frac{2^{n+3}}{3^{n+1}z^2}\left(1-\frac{2}{3}\sin^2\frac{z}{2}\right)^{-1}\sin^{2n+4}\frac{z}{2}\right) +$$

$$+\left(\frac{\lambda}{\kappa}\right)^2\frac{\sin z}{p^2 z}\left(2\frac{\sin z}{z}\left(1-\frac{2}{3}\sin^2\frac{z}{2}\right)^{-1} - 2 - \frac{2^{n+1}\sin z}{3^{n+1}z}\left(1-\frac{2}{3}\sin^2\frac{z}{2}\right)^{-1}\sin^{2n+2}\frac{z}{2}\right).$$

Moreover, $g_n(z) \geq f_n(z)$ for arbitrary $z \neq 0$, where $f_n(z)$ is defined as in Lemma 1.

**Proof.** Let $A_n$ and $B_n$ are defined in the same way as in the proof of Lemma 1 above. Then

$$\omega_G = \left(-\frac{4\kappa}{h^2}\sin^2\left(\frac{ph}{2}\right) - i\frac{\lambda}{h}\sin(ph)\right)\sum_{n=0}^{+\infty}\frac{2^n}{3^n}\sin^{2n}\left(\frac{ph}{2}\right) = \omega_n - (A_{n+1} + A_{n+2} + \ldots) - i(B_{n+1} + B_{n+2} + \ldots)$$

(the infinite geometric progression with the common ratio $(2/3)\sin^2(ph/2)$), and thereby

$$|\omega_G - \overline{\omega}|^2 - |\omega_n - \overline{\omega}|^2 = (A_{n+1} + A_{n+2} + \ldots)((A_{n+1} + A_{n+2} + \ldots) - 2\operatorname{Re}(\omega_n - \overline{\omega})) + (B_{n+1} + B_{n+2} + \ldots) \times$$

$$\times((B_{n+1} + B_{n+2} + \ldots) - 2\operatorname{Im}(\omega_n - \overline{\omega})) = \frac{2^{n+3}\kappa}{3^{n+1}h^2}\left(1 - \frac{2}{3}\sin^2\frac{ph}{2}\right)^{-1}\sin^{2n+4}\left(\frac{ph}{2}\right) \times$$

$$\times\left((A_{n+1} + A_{n+2} + \ldots) + 2\left(\frac{4\kappa}{h^2}\sin^2\frac{ph}{2} + \frac{8\kappa}{3h^2}\sin^4\frac{ph}{2} + \ldots + \frac{2^{n+2}\kappa}{3^n h^2}\sin^{2n+2}\frac{ph}{2} - \kappa p^2\right)\right) + \frac{2^{n+1}\lambda}{3^{n+1}h}\sin(ph) \times$$

$$\times\left((B_{n+1} + B_{n+2} + \ldots) + 2\frac{\lambda\sin(ph)}{h}\left(1 + \frac{2}{3}\sin^2\frac{ph}{2} + \ldots + \frac{2^n}{3^n}\sin^{2n}\frac{ph}{2}\right) - 2\lambda p\right)\left(1 - \frac{2}{3}\sin^2\frac{ph}{2}\right)^{-1}\sin^{2n+2}\frac{ph}{2} =$$

$$= \frac{2^{n+1}\kappa^2 p^4}{3^{n+1}}\sin^{2n+2}\left(\frac{z}{2}\right)\left(1 - \frac{2}{3}\sin^2\frac{z}{2}\right)^{-1}g_n(z),$$

$$g_n(z) = \frac{4}{z^2}\sin^2\left(\frac{z}{2}\right)\left(\frac{4}{z^2}\sin^2\frac{z}{2}\left(\sum_{m=n+1}^{+\infty}\frac{2^m}{3^m}\sin^{2m}\frac{z}{2}\right) + \frac{8}{z^2}\sin^2\frac{z}{2}\left(\sum_{m=0}^{n}\frac{2^m}{3^m}\sin^{2m}\frac{z}{2}\right) - 2\right) +$$

$$+\left(\frac{\lambda}{\kappa}\right)^2\frac{\sin z}{p^2 z}\left(\frac{\sin z}{z}\left(\sum_{m=n+1}^{+\infty}\frac{2^m}{3^m}\sin^{2m}\frac{z}{2}\right) + 2\frac{\sin z}{z}\left(\sum_{m=0}^{n}\frac{2^m}{3^m}\sin^{2m}\frac{z}{2}\right) - 2\right) =$$

$$= \frac{4}{z^2}\sin^2\left(\frac{z}{2}\right)\left(\frac{8}{z^2}\sin^2\frac{z}{2}\left(1-\frac{2}{3}\sin^2\frac{z}{2}\right)^{-1} - 2 - \frac{2^{n+3}}{3^{n+1}z^2}\left(1-\frac{2}{3}\sin^2\frac{z}{2}\right)^{-1}\sin^{2n+4}\frac{z}{2}\right) +$$

$$+\left(\frac{\lambda}{\kappa}\right)^2\frac{\sin z}{p^2 z}\left(2\frac{\sin z}{z}\left(1-\frac{2}{3}\sin^2\frac{z}{2}\right)^{-1} - 2 - \frac{2^{n+1}\sin z}{3^{n+1}z}\left(1-\frac{2}{3}\sin^2\frac{z}{2}\right)^{-1}\sin^{2n+2}\frac{z}{2}\right).$$

$$g_n(z) - f_n(z) = \frac{4}{z^2}\sin^2\frac{z}{2}\left(\frac{4}{z^2}\sin^2\frac{z}{2}\left(\sum_{m=n+2}^{+\infty}\frac{2^m}{3^m}\sin^{2m}\frac{z}{2}\right)\right) + \frac{\lambda^2\sin z}{\kappa^2 p^2 z}\left(\frac{\sin z}{z}\left(\sum_{m=n+2}^{+\infty}\frac{2^m}{3^m}\sin^{2m}\frac{z}{2}\right)\right) \geq 0. \blacksquare$$

**Remark 2.** Let $\kappa = 0$, $n \geq 1$ and $\lambda \neq 0$. Then we have the following equality (where $z \equiv ph \neq 0$):

$$|\omega_G - \overline{\omega}|^2 - |\omega_n - \overline{\omega}|^2 = \left(\frac{2^{n+1}\lambda^2 p^2}{3^{n+1}}\left(1-\frac{2}{3}\sin^2\frac{z}{2}\right)^{-1}\sin^{2n+2}\frac{z}{2}\right)\tilde{g}_n(z), \text{ where } \tilde{g}_n(z) \text{ is defined by}$$

$$\tilde{g}_n(z) = \frac{\sin z}{z}\left(\left(1-\frac{2}{3}\sin^2\frac{z}{2}\right)^{-1}\frac{\sin z}{z}\left(2 - \frac{2^{n+1}}{3^{n+1}}\sin^{2n+2}\frac{z}{2}\right) - 2\right).$$

Moreover, $\tilde{g}_n(z) \leq 2 \dfrac{\sin z}{z} \left( \dfrac{\sin z}{z} \left( 1 - \dfrac{2}{3} \sin^2 \dfrac{z}{2} \right)^{-1} - 1 \right)$ for arbitrary $n \geq 1$ and $z \neq 0$.

**Proof.** Following the proof of Lemma 2, it is easy to see that in this case we have $|\omega_G - \overline{\omega}|^2 - |\omega_n - \overline{\omega}|^2 =$

$$= \left( \dfrac{\sin^2 z}{z^2} \left( \sum_{m=n+1}^{+\infty} \dfrac{2^m}{3^m} \sin^{2m} \dfrac{z}{2} \right) + 2 \dfrac{\sin^2 z}{z^2} \left( \sum_{m=0}^{n} \dfrac{2^m}{3^m} \sin^{2m} \dfrac{z}{2} \right) - 2 \dfrac{\sin z}{z} \right) \dfrac{2^{n+1} \lambda^2 p^2}{3^{n+1}} \left( 1 - \dfrac{2}{3} \sin^2 \dfrac{z}{2} \right)^{-1} \sin^{2n+2} \dfrac{z}{2},$$

which is equivalent to the above representation. Next, we immediately have the following inequality:

$$2 \dfrac{\sin z}{z} \left( \dfrac{\sin z}{z} \left( 1 - \dfrac{2}{3} \sin^2 \dfrac{z}{2} \right)^{-1} - 1 \right) - \tilde{g}_n(z) = \dfrac{\sin^2 z}{z^2} \left( 1 - \dfrac{2}{3} \sin^2 \dfrac{z}{2} \right)^{-1} \dfrac{2^{n+1}}{3^{n+1}} \sin^{2n+2} \dfrac{z}{2} \geq 0. \blacksquare$$

For illustrative purposes, the functions $\tilde{f}_i(z)$ and $\tilde{g}_i(z)$ (for $1 \leq i \leq 4$ and $0 < z \leq \pi$) are plotted in Figure 1a. The functions $f_i(z)$, $g_i(z)$ ($1 \leq i \leq 4$, $0 < z \leq \pi$) depending on $\mu \equiv \lambda/(\kappa p)$ for $\mu = 0$, $\mu = 1$, $\mu = 2$ and $\mu = \lambda/(\kappa p) \approx 10.6$ (where $\lambda = 1$, $\kappa = 10^{-2}$, $p = 3\pi$) are plotted in Figure 1b, Figure 2a, Figure 2b and 3a, respectively. Figure 3b shows a close-up view of Figure 3a for $0 < z \leq 2z_0$ ($z_0 \approx 0.195$ in this case, where $z_0$ is defined as in Lemma 1).

The number $z_0$ plays an important theoretical role – it shows that for arbitrary parameters $\kappa > 0$, $\lambda$ and $p \neq 0$, there is always a neighbourhood of zero where the functions $f_n(z)$ and $g_n(z)$ ($n \geq 1$) are guaranteed to be positive (respectively, then $|\omega_{n+1} - \overline{\omega}| > |\omega_n - \overline{\omega}|$ and $|\omega_G - \overline{\omega}| > |\omega_n - \overline{\omega}|$ – see details in Propositions 3 and 5 below). Later we will study the properties of $z_0$ in more detail (see Remark 3 below). In particular, we will show that the number $z_0$ can be considered as an approximation of the smallest positive root of $f_1(z)$ from the left side and they rapidly approach each other under certain circumstances. Note that $z_0$ was chosen to be the root of the first suitable polynomial such that the graph of $f_1(z)$ is guaranteed to lie above the graph of this polynomial in a zero vicinity for arbitrary admissible $\kappa$, $\lambda$ and $p$; thus applying the technique and power expansions used in the proof of Lemma 1, one can find higher-degree polynomials whose roots approximate the root of the function $f_1(z)$ even more accurately, but there is no special theoretical meaning in this, since in general high-degree equations cannot be solved algebraically (owing to the well-known Abel-Ruffini theorem).

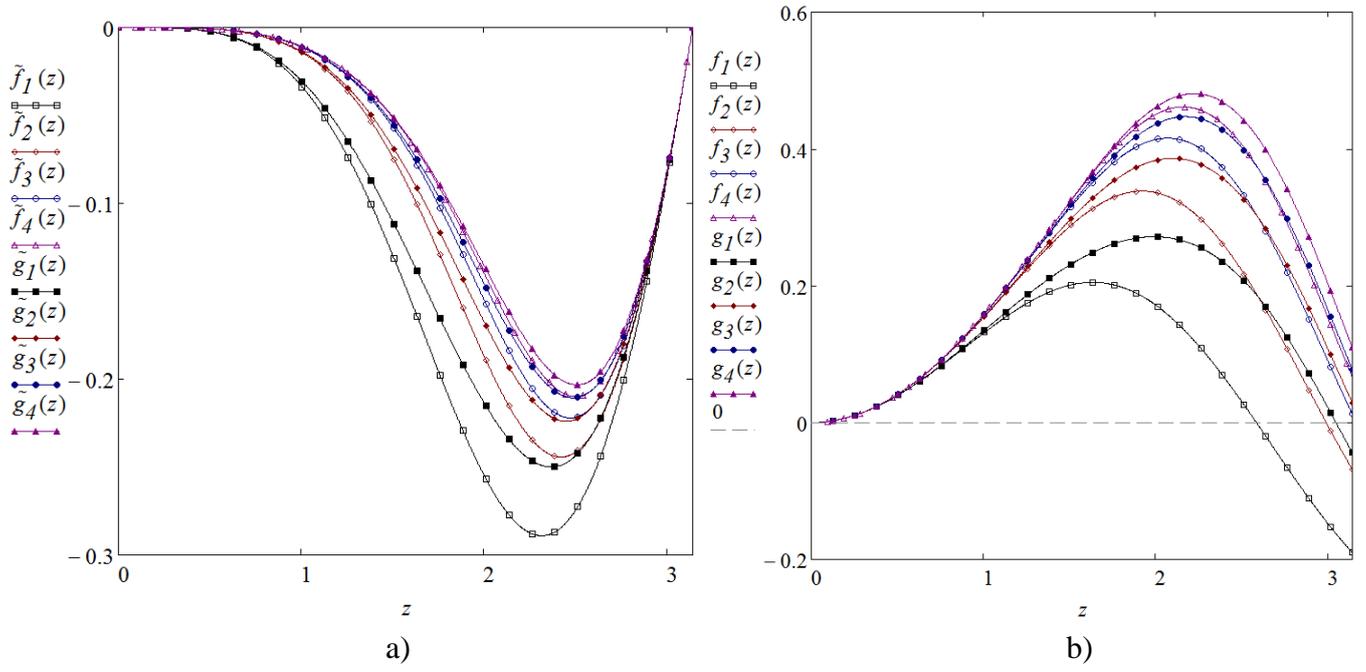

Figure 1. a) $\widetilde{f}_i(z)$ and $\widetilde{g}_i(z)$, $1 \leq i \leq 4$, $0 < z \leq \pi$; b) $f_i(z)$ and $g_i(z)$, $\mu = 0$, $1 \leq i \leq 4$, $0 < z \leq \pi$.

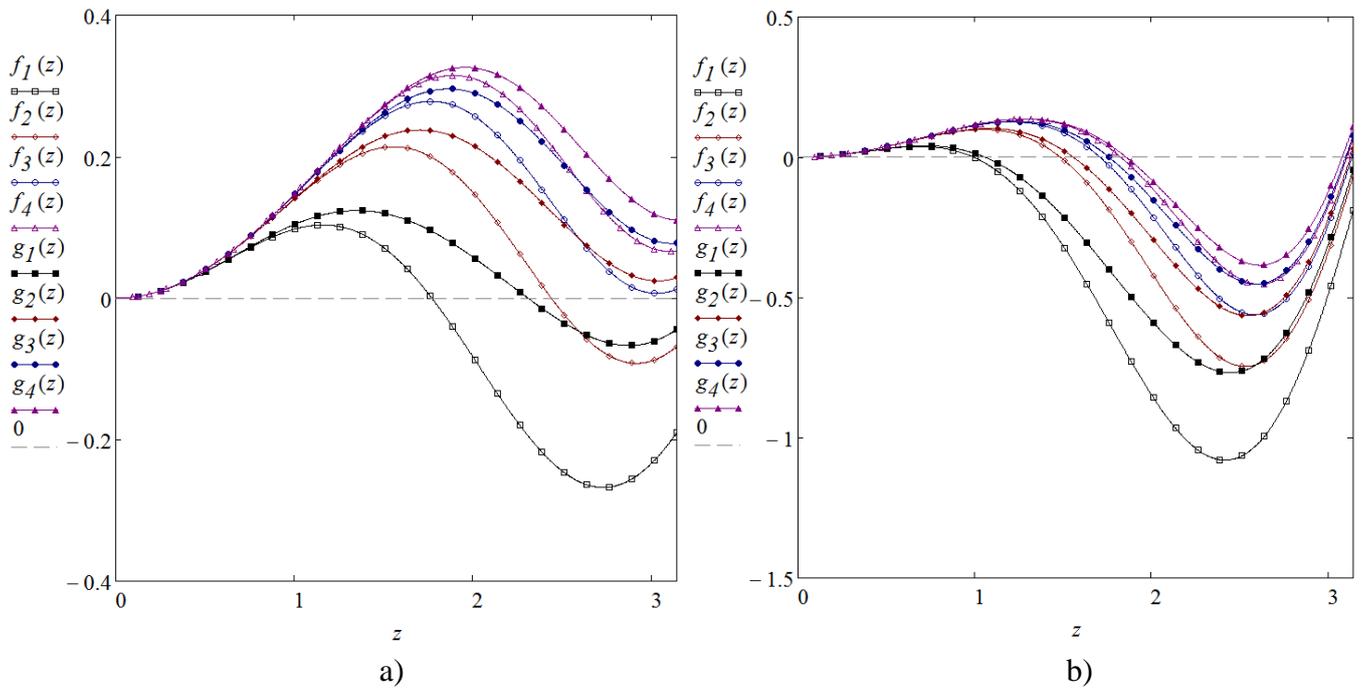

Figure 2. a) $f_i(z)$ and $g_i(z)$, $\mu = 1$, $1 \leq i \leq 4$, $0 < z \leq \pi$; b) $f_i(z)$, $g_i(z)$, $\mu = 2$, $1 \leq i \leq 4$, $0 < z \leq \pi$.

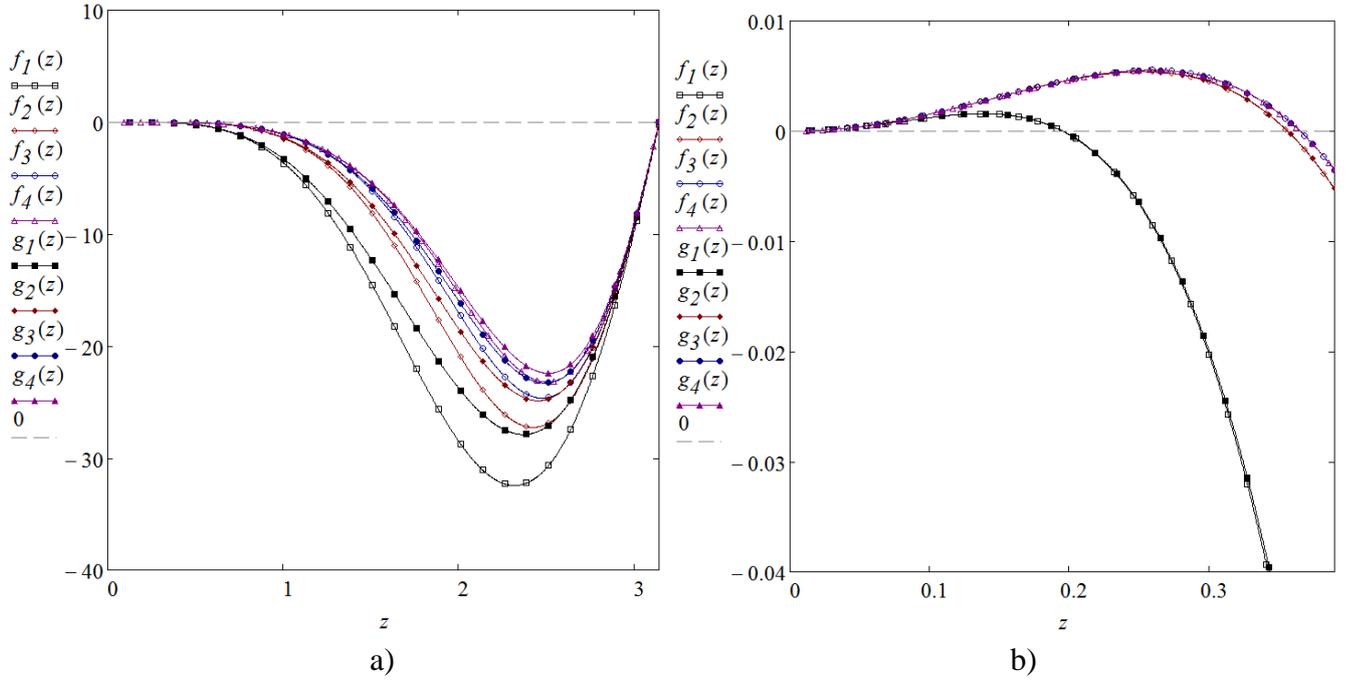

*Figure 3.* $f_i(z)$, $g_i(z)$, $\mu = 10.6$, $1 \leq i \leq 4$: *a)* $0 < z \leq \pi$; *b) a close-up view on* $0 < z \leq 2z_0$.

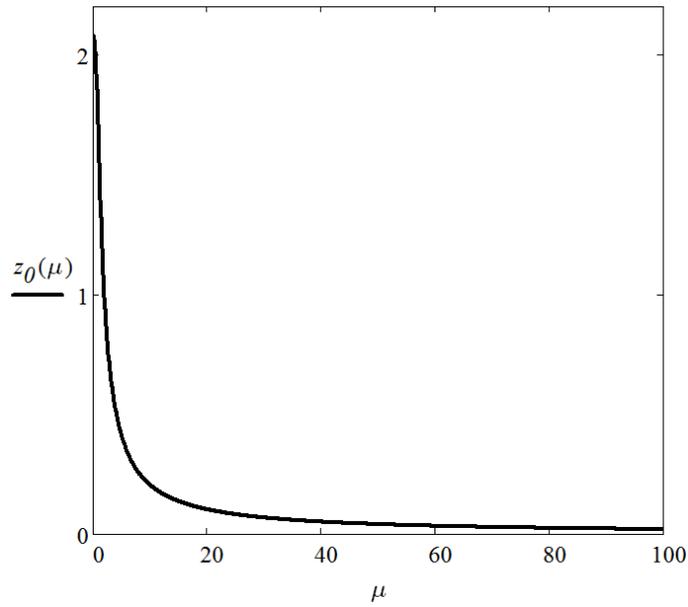

*Figure 4.* $z_0 = z_0(\mu)$ *where* $\mu = \lambda/(\kappa p)$.

**Proposition 3.** Let $\kappa > 0$, $n \geq 1$, and let $\lambda$ be arbitrary. Then for arbitrary $p \neq 0$ the inequality $|\omega_n - \bar{\omega}| < |\omega_{n+1} - \bar{\omega}|$ is true for all $h$ satisfying $h \leq z_0/|p|$ independently of $n$, where $z_0$ is defined as in Lemma 1.

**Proof.** The proof follows directly from Lemma 1. Indeed, owing to Lemma 1 we have $z_0 < \pi$ and

$$|\omega_{n+1} - \bar{\omega}|^2 - |\omega_n - \bar{\omega}|^2 \geq \left( \frac{2^{n+1} \kappa^2 p^4}{3^{n+1}} \sin^{2n+2}\left(\frac{z}{2}\right) \right) \cdot f_1(z) > 0 \text{ for all } z = ph \text{ satisfying } 0 < |z| \leq z_0. \blacksquare$$

**Proposition 4.** Let $\kappa = 0$, $n \geq 1$, and $\lambda \neq 0$. Then for arbitrary $p \neq 0$ the inequality $|\omega_{n+1} - \overline{\omega}| < |\omega_n - \overline{\omega}|$ is true for all $h$ satisfying $h < \pi/|p|$ independently of $n$.

**Proof.** Using the representation for $|\omega_{n+1} - \overline{\omega}|^2 - |\omega_n - \overline{\omega}|^2$ from Remark 1, we obtain

$$|\omega_{n+1} - \overline{\omega}|^2 - |\omega_n - \overline{\omega}|^2 \leq 2 \frac{\sin z}{z} \left( \frac{\sin z}{z} \left( 1 - \frac{2}{3} \sin^2 \frac{z}{2} \right)^{-1} - 1 \right) \frac{2^{n+1} \lambda^2 p^2}{3^{n+1}} \sin^{2n+2} \frac{z}{2} < 0, \text{ since the function}$$

$\frac{\sin z}{z}$ is strictly positive and the function $\frac{\sin z}{z} \left( 1 - \frac{2}{3} \sin^2 \frac{z}{2} \right)^{-1} - 1$ is strictly negative when $0 < |z| < \pi$. ∎

**Proposition 5.** Let $\kappa > 0$, $n \geq 1$, and let $\lambda$ be arbitrary. Then for arbitrary $p \neq 0$ the inequality $|\omega_n - \overline{\omega}| < |\omega_G - \overline{\omega}|$ is true for all $h$ satisfying $h \leq z_0/|p|$ independently of $n$, where $z_0$ is defined as in Lemma 1.

**Proof.** The proof follows directly from Lemma 1 and Lemma 2. Indeed, owing to Lemma 2 we have

$$|\omega_G - \overline{\omega}|^2 - |\omega_n - \overline{\omega}|^2 = \left( \frac{2^{n+1} \kappa^2 p^4}{3^{n+1}} \left( 1 - \frac{2}{3} \sin^2 \frac{z}{2} \right)^{-1} \sin^{2n+2} \frac{z}{2} \right) g_n(z), \text{ where } g_n(z) \geq f_n(z). \text{ Since}$$

$f_n(z) \geq f_1(z)$, $f_1(z) > 0$ when $0 < |z| \leq z_0$ and $z_0 < \pi$ (see Lemma 1), the inequality $|\omega_n - \overline{\omega}| < |\omega_G - \overline{\omega}|$ now follows immediately. ∎

**Proposition 6.** Let $\kappa = 0$, $n \geq 1$, and $\lambda \neq 0$. Then for arbitrary $p \neq 0$ the inequality $|\omega_G - \overline{\omega}| < |\omega_n - \overline{\omega}|$ is true for all $h$ satisfying $h < \pi/|p|$ independently of $n$.

**Proof.** Using the representation for $|\omega_G - \overline{\omega}|^2 - |\omega_n - \overline{\omega}|^2$ from Remark 2, we immediately obtain

$$|\omega_G - \overline{\omega}|^2 - |\omega_n - \overline{\omega}|^2 \leq 2 \frac{\sin z}{z} \left( \frac{\sin z}{z} \left( 1 - \frac{2}{3} \sin^2 \frac{z}{2} \right)^{-1} - 1 \right) \frac{2^{n+1} \lambda^2 p^2}{3^{n+1}} \left( 1 - \frac{2}{3} \sin^2 \frac{z}{2} \right)^{-1} \sin^{2n+2} \frac{z}{2} < 0, \text{ since}$$

$\frac{\sin z}{z}$ is strictly positive and the function $\frac{\sin z}{z} \left( 1 - \frac{2}{3} \sin^2 \frac{z}{2} \right)^{-1} - 1$ is strictly negative when $0 < |z| < \pi$. ∎

**Remark 3.** In accordance with the above results (see Lemmas 1-2 and Remarks 1-2), the signs of $|\omega_{n+1} - \overline{\omega}| - |\omega_n - \overline{\omega}|$ and $|\omega_G - \overline{\omega}| - |\omega_n - \overline{\omega}|$ are completely determined by the signs of the functions $f_n(z)$ and $g_n(z)$ ($\tilde{f}_n(z)$ and $\tilde{g}_n(z)$ for the case $\kappa = 0$), respectively. Recall that the convection-diffusion ratio is usually characterized by the so-called Peclet number [1, 2, 5–7], which is an important physical parameter that shows the rate of the prevalence of convection over diffusion. Lemma 1 states that for an arbitrary finite Peclet number $Pe \equiv |\lambda|/\kappa < \infty$ and arbitrary $p \neq 0$ there exists $z_0 < \pi$ such that the inequality $0 < |z| \leq z_0$ gives *a sufficient condition* for positivity of the functions $f_n(z)$, $n \geq 1$. Respectively, there are values of $h$ (namely, $h \leq z_0/|p|$) such that $|\omega_n - \overline{\omega}| < |\omega_{n+1} - \overline{\omega}|$ and $|\omega_n - \overline{\omega}| < |\omega_G - \overline{\omega}|$ independently of $n$, see Propositions 3 and 5 for details. But note that the region for such $h$ narrows and decreases with increasing the number $Pe$ (i.e., with growth in the dominance of convection over diffusion) – this follows from the facts that:

1) *The number $z_0$ tends monotonically to zero from the right as the number Pe passes from 0 to $\infty$. Moreover, $z_0 = O((Pe)^{-1})$ as $Pe \to \infty$ and the convergence rate $-1$ is exact here;*

2) $\forall Pe < \infty$ $\forall p \neq 0$ $\exists \tilde{z} \in (z_0, \pi)$ that $f_1(\tilde{z}) = 0$ and $f_1(z) > 0$ for all $0 < |z| < \tilde{z}$, and $\lim_{Pe \to \infty} \tilde{z} = 0$. Moreover, $\tilde{z} = O((Pe)^{-1})$ (the convergence rate $-1$ is exact here) and $\tilde{z} - z_0 = O((Pe)^{-3})$ as $Pe \to \infty$. The number $z_0$ depending on $\mu \equiv \lambda/(\kappa p)$ is plotted in Figure 4 (note also that $|\mu| = Pe/|p|$ and the function $z_0 = z_0(\mu)$ is even). The second item means that the number $\tilde{z}$ is the smallest positive root of the function $f_1(z)$, passing through which the function changes its sign from plus to minus, and that this root tends to zero from the right side as $Pe \to \infty$. Finally, for pure transport problems with $\kappa = 0$ (formally this is the case of $Pe = \infty$ or the absolute dominance of convection over diffusion) we obtain the inverse relations – that is, $|\omega_{n+1} - \overline{\omega}| < |\omega_n - \overline{\omega}|$ and $|\omega_G - \overline{\omega}| < |\omega_n - \overline{\omega}|$ independently of $n$ and for all $h < \pi/|p|$, see Propositions 4 and 6 for details.

**Proof.** All statements of this remark are obvious, except for items 1 and 2, which we shall now prove. A direct calculation shows that $z_0 \xrightarrow[Pe \to \infty]{} 0$ monotonically, where $z_0$ is defined as in Lemma 1, and $\exists \lim_{Pe \to \infty} z_0/(Pe)^{-1} = \sqrt{210}|p|/7 \neq 0$. Redefining the function $f_1(z)$ at $z = 0$ by its limit, $f_1(0) \equiv \lim_{z \to 0} f_1(z) = 0$, we obtain that $f_1(z)$ is continuous on $[0, \pi]$. Since $f_1(z_0) > 0$ (see the proof of Lemma 1) and $f_1(\pi) = 8(68 - 9\pi^2)/9\pi^4 < 0$ for arbitrary $Pe$ and $p$, then $\exists \tilde{z} \in (z_0, \pi)$ (due to the continuity of $f_1(z)$ [26]) that $f_1(\tilde{z}) = 0$ and $f_1(z) > 0$ when $z \in (z_0, \tilde{z})$. Let us now prove that a root of $f_1(z)$ with such properties (i.e., the smallest positive root of $f_1(z)$) tends to zero as $Pe \to \infty$. The basic idea of this proof lies in the following: we find some point $\xi$ such that $z_0 < \xi$, $f_1(\xi) < 0$ and $\lim_{Pe \to \infty} \xi = 0$. Then by virtue of continuity of $f_1(z)$ and known relations $f_1(z_0) > 0$ and $\lim_{Pe \to \infty} z_0 = 0$ we immediately obtain that $\tilde{z} \in (z_0, \xi)$ and $\lim_{Pe \to \infty} \tilde{z} = 0$, that is the desired result. Using the expansion for $f_1(z)$ and the alternating series test (see the proof of Lemma 1) we obtain the inequality $P_1(z) < f_1(z) < P_2(z)$ if $0 < |z| < z_*$, where $z_*$ and $c_m$ are defined as in the proof of Lemma 1, and the polynomials $P_1$, $P_2$ are defined by $P_1(z) \equiv z^2(5040 - 12z^2(98\mu^2 + 91) + z^4(156\mu^2 - 17))/30240$ and $P_2(z) \equiv P_1(z) + c_4 z^8 = P_1(z) + z^8(79\mu^2 + 80)/100800$, $\mu \equiv \lambda/(\kappa p)$. Let us then introduce the polynomial $\tilde{P}_2(z) \equiv P_1(z) + z^6(79\mu^2 + 80)/100800$ and denote by $\psi$ the smallest positive root of $\tilde{P}_2(z)$, namely

$$\psi = \sqrt{\frac{5880\mu^2 + 5460 - \sqrt{(5880\mu^2 + 5460)^2 - 50400(1797\mu^2 + 70)}}{1797\mu^2 + 70}} =$$

$$= \sqrt{\frac{5880\mu^2 + 5460 - 60\sqrt{7}\sqrt{7(14\mu^2 - 523/98)^2 + 1157467/1372}}{1797\mu^2 + 70}}.$$

It is easy to see that $\psi$ exists for all $\mu$, and $\lim_{Pe \to \infty} \psi = 0$ (note that $|\mu| = Pe/|p|$). Since we now investigate the situation of $Pe \to \infty$ (respectively, $|\mu| = Pe/|p| \to \infty$ as well), we can assume without loss of generality that $\mu^2 > 39550/9963$. Then a direct arithmetic calculation shows that this condition ensures the inequality $z_0 < \psi < 1 < z_*$, thereby providing $f_1(\psi) < P_2(\psi) < \tilde{P}_2(\psi) = 0$. The last result shows that one can take $\psi$ as the required number $\xi$. Now the asymptotical relation $\tilde{z} - z_0 = O((Pe)^{-3})$ follows from the facts that $(\tilde{z} - z_0)/(Pe)^{-3} < (\xi - z_0)/(Pe)^{-3}$ and the expression $(\xi - z_0)/(Pe)^{-3}$ has a finite limit as $Pe \to \infty$ (a direct calculation shows that $\lim_{Pe \to \infty} (\xi - z_0)/(Pe)^{-3} = 237\sqrt{210}|p|^3/38416 \neq 0$). Similarly, $\exists \lim_{Pe \to \infty} \tilde{z}/(Pe)^{-1} = \sqrt{210}|p|/7 \neq 0$ since $\lim_{Pe \to \infty} \xi/(Pe)^{-1} = \lim_{Pe \to \infty} z_0/(Pe)^{-1} = \sqrt{210}|p|/7$ and $z_0 < \tilde{z} < \xi$, which implies the exact asymptotic relation $\tilde{z} = O((Pe)^{-1})$ as $Pe \to \infty$. ∎

Let us denote the number $\omega$ for the standard (non-corrected) lumped formulation by $\omega_L$. Then, substituting the ansatz $a_k(t) = e^{\omega_L t} e^{ikhp}$ into (3), we obtain the following expression for $\omega_L$:

$$\omega_L = -\frac{4\kappa}{h^2}\sin^2\left(\frac{ph}{2}\right) - i\frac{\lambda}{h}\sin(ph).$$

It easy to see that $\lim_{h \to 0} \omega_L = \overline{\omega}$ and that $\omega_L = \omega_0$ if we formally put $n = 0$ in the expression for $\omega_n$ (see Proposition 2). Let us note that the case of using the first correction (i.e. the case with $n = 1$ in Definition 1) is the most frequently used case of application of the Guermond-Pasquetti technique in literature (see [12–17]), mainly due to its greatest simplicity and observations in [11] that the first correction is usually sufficient to compensate the dominating dispersive effects of mass lumping in pure transport problems (in particular, it is rigorously proved in [11] that the first correction eliminates the leading terms in the consistency error of lumped scheme in the 1D pure transport case). The following proposition justifies this via Fourier analysis (note that we also consider the case with the presence of diffusion).

**Proposition 7.** Let parameters $\kappa$ and $\lambda$ be arbitrary, $\kappa + |\lambda| > 0$. Then for arbitrary $p \neq 0$ the inequality $|\omega_1 - \overline{\omega}| < |\omega_L - \overline{\omega}|$ is true for all $h$ satisfying $h < \pi/|p|$.

**Proof.** Let $A_n$ and $B_n$ are defined as in the proof of Lemma 1. Following the proof of Lemma 1 we have

$$|\omega_1 - \overline{\omega}|^2 - |\omega_0 - \overline{\omega}|^2 = |(\omega_0 - A_1 - iB_1) - \overline{\omega}|^2 - |\omega_0 - \overline{\omega}|^2 = A_1(A_1 - 2\operatorname{Re}(\omega_0 - \overline{\omega})) + B_1(B_1 - 2\operatorname{Im}(\omega_0 - \overline{\omega})) =$$

$$= \left(\frac{4p^2}{3}\sin^2\left(\frac{z}{2}\right)\right)\left(\kappa^2 p^2 \frac{4}{z^2}\left(\frac{4}{z^2}\sin^2\frac{z}{2} + \frac{4}{3z^2}\sin^4\frac{z}{2} - 1\right)\sin^2\frac{z}{2} + \lambda^2 \frac{\sin z}{z}\left(\frac{\sin z}{z} + \frac{\sin z}{3z}\sin^2\frac{z}{2} - 1\right)\right) < 0$$

for all $|z| = |ph| \in (0, \pi)$ since the functions $\left(\frac{4}{z^2}\sin^2\frac{z}{2} + \frac{4}{3z^2}\sin^4\frac{z}{2} - 1\right)$ and $\left(\frac{\sin z}{z} + \frac{\sin z}{3z}\sin^2\frac{z}{2} - 1\right)$ are strictly negative if $0 < |z| < \pi$. ∎

**Remark 4.** Note that the expression $|\omega_1 - \overline{\omega}|^2 - |\omega_0 - \overline{\omega}|^2$ derived in the proof of Proposition 7 coincides with the representations for $|\omega_{n+1} - \overline{\omega}|^2 - |\omega_n - \overline{\omega}|^2$ derived in Lemma 1 (for $\kappa > 0$) and Remark 1 (for $\kappa = 0$) if we formally put $n = 0$ there. However, as Proposition 7 shows, the behaviour of $f_0(z)$ is completely opposite to the behaviour of the functions $f_n(z)$ ($n \geq 1$) in the sense that the function $f_0(z)$ is always negative for all $0 < |z| < \pi$. This shows the important role of the additional terms

$$-\frac{8\kappa}{3h^2}\sin^4\left(\frac{ph}{2}\right) - \frac{16\kappa}{9h^2}\sin^6\left(\frac{ph}{2}\right) - \ldots - \frac{2^{n+2}\kappa}{3^n h^2}\sin^{2n+2}\left(\frac{ph}{2}\right) \text{ in } \omega_n \text{ for the existence the region (namely,}$$

$0 < |z| \leq z_0$ − see Lemma 1) where the functions $f_n(z)$ ($n \geq 1$) are a-priori positive.

**Proposition 8.** Let $n \geq 1$, and parameters $\kappa$ and $\lambda$ be arbitrary, $\kappa + |\lambda| > 0$. Then for arbitrary $p \neq 0$ the inequality $|\omega_G - \omega_{n+1}| < |\omega_G - \omega_n|$ is true for all $h$ satisfying $h < \pi/|p|$ independently of $n$.

**Proof.** Let $A_n$ and $B_n$ are defined in the same way as in the proof of Lemma 1 above. Using the representation $\omega_G = \omega_n - \sum_{m=n+1}^{+\infty}(A_m + iB_m)$ (see the proof of Lemma 2), we immediately obtain

$$|\omega_G - \omega_n|^2 - |\omega_G - \omega_{n+1}|^2 = \left(\sum_{m=n+1}^{+\infty} A_m\right)^2 + \left(\sum_{m=n+1}^{+\infty} B_m\right)^2 - \left(\sum_{m=n+2}^{+\infty} A_m\right)^2 - \left(\sum_{m=n+2}^{+\infty} B_m\right)^2 = (A_{n+1})^2 + (B_{n+1})^2 +$$

$$+ 2A_{n+1} \cdot \sum_{m=n+2}^{+\infty} A_m + 2B_{n+1} \cdot \sum_{m=n+2}^{+\infty} B_m = \left(\frac{2^{n+3}\kappa}{3^{n+1}h^2}\sin^{2n+4}\left(\frac{ph}{2}\right)\right)^2 + \left(\frac{2^{n+1}\lambda}{3^{n+1}h}\sin(ph)\sin^{2n+2}\left(\frac{ph}{2}\right)\right)^2 +$$

$$+ \frac{2^{n+4}\kappa^2}{3^{n+1}h^4}\sin^{2n+4}\left(\frac{ph}{2}\right)\left(\sum_{m=n+2}^{+\infty}\frac{2^{m+2}}{3^m}\sin^{2m+2}\left(\frac{ph}{2}\right)\right) + \frac{2^{n+2}\lambda^2\sin^2(ph)}{3^{n+1}h^2}\sin^{2n+2}\left(\frac{ph}{2}\right)\left(\sum_{m=n+2}^{+\infty}\frac{2^m}{3^m}\sin^{2m}\left(\frac{ph}{2}\right)\right),$$

which is strictly positive if $0 < |ph| < \pi$. ∎

In the following proposition we use the standard concept of the (local) equivalence of functions in mathematical analysis (e.g., see [26]): namely, $f$ is equivalent to $g$ as $h \to 0$ (we shall write $f \sim g$ as $h \to 0$) if and only if $\lim_{h \to 0}(f/g) = 1$. This means that $f$ behaves asymptotically like $g$ as $h \to 0$ [26].

**Proposition 9** (asymptotical formulas). 1) Let $\kappa > 0$, $n \geq 1$, $p \neq 0$ and let $\lambda$ be arbitrary, then

$$|\omega_{n+1} - \overline{\omega}|^2 - |\omega_n - \overline{\omega}|^2 \sim \frac{\kappa^2 p^{2n+8}}{6^{n+2}}h^{2n+4} \quad \text{and} \quad |\omega_G - \overline{\omega}|^2 - |\omega_n - \overline{\omega}|^2 \sim \frac{\kappa^2 p^{2n+8}}{6^{n+2}}h^{2n+4} \quad \text{as } h \to 0.$$

2) Let $\kappa = 0$, $n \geq 1$, $p \neq 0$ and $\lambda \neq 0$, then $|\omega_{n+1} - \overline{\omega}|^2 - |\omega_n - \overline{\omega}|^2 \sim -\lambda^2 p^{2n+8} a_n h^{2n+6}$ and $|\omega_G - \overline{\omega}|^2 - |\omega_n - \overline{\omega}|^2 \sim -\lambda^2 p^{2n+8} a_n h^{2n+6}$ as $h \to 0$, where $a_1 = 7/6480$, $a_n = 1/(15 \cdot 6^{n+2})$ ($n \geq 2$).

**Proof.** A direct calculation of limits with the usage of representations given in Lemmas 1-2 and Remarks 1-2 immediately yields the corresponding results. ∎

**Remark 5.** Proposition 9 implies that $|\omega_{n+1} - \overline{\omega}|^2 - |\omega_n - \overline{\omega}|^2 = O(h^{2n+m})$ and $|\omega_G - \overline{\omega}|^2 - |\omega_n - \overline{\omega}|^2 = O(h^{2n+m})$ as $h \to 0$, where $m = 4$ if $\kappa > 0$ and $m = 6$ if $\kappa = 0$, and that the order of convergence $2n + m$ is exact (unimprovable). Note also that the expressions for $|\omega_{n+1} - \overline{\omega}|^2 - |\omega_n - \overline{\omega}|^2$ and $|\omega_G - \overline{\omega}|^2 - |\omega_n - \overline{\omega}|^2$ (obtained in Lemmas 1-2 and Remarks 1-2) are even with respect to $z$, so their power expansions contain only even powers of $h$ (respectively, Proposition 9 provides the first terms of the corresponding expansions).

**Remark 6.** Owing to the complexity of direct investigation of exponentials in harmonics, the study of the interrelations between $\omega_{\text{numerical}}$ and $\omega_{\text{exact}} = \overline{\omega}$ (in particular, the estimation of the distance between them or various ratios between their absolute values, real/imaginary parts, etc.) is a conventional way in Fourier analysis of numerical approximations [1, 2, 18–20, 22, 23, 27–29]. Using asymptotical formulas for exponentials, and the equality $|e^{ikhp}| = 1$ ($khp$ is real), it is easy to show (see [22, 23, 29]) that $\|u_{\text{numerical}} - u_{\text{exact}}\|_\infty = \max_{x_k}|u_{\text{numerical}}(t, x_k) - u_{\text{exact}}(t, x_k)| = \max_k|e^{t\omega_{\text{numerical}}}e^{ikhp} - e^{t\omega_{\text{exact}}}e^{ikhp}| =$

$= \left| e^{t\omega_{\text{numerical}}} - e^{t\omega_{\text{exact}}} \right| \sim \left| e^{t\omega_{\text{exact}}} \right| \left| \omega_{\text{numerical}} - \omega_{\text{exact}} \right| t$ as $h \to 0$, which emphasizes the above-mentioned important role of the distance $\left| \omega_{\text{numerical}} - \omega_{\text{exact}} \right|$ in Fourier analysis.

**Remark 7.** Let us denote $\tilde{R}_n \equiv \left\| e^{t\omega_{n+1}} e^{ikhp} - e^{\overline{\omega}t} e^{ikhp} \right\|_\infty^2 - \left\| e^{t\omega_n} e^{ikhp} - e^{\overline{\omega}t} e^{ikhp} \right\|_\infty^2$ and $\tilde{G}_n \equiv \left\| e^{t\omega_G} e^{ikhp} - e^{\overline{\omega}t} e^{ikhp} \right\|_\infty^2 - \left\| e^{t\omega_n} e^{ikhp} - e^{\overline{\omega}t} e^{ikhp} \right\|_\infty^2$. For the convenience of verifying the numerical results of Section 3 below, let us now show that

$$\tilde{R}_n = O(h^{2n+m}) \text{ and } \tilde{G}_n = O(h^{2n+m}) \text{ as } h \to 0 \text{ for arbitrary fixed } t > 0,$$

where $m = 4$ if $\kappa > 0$ and $m = 6$ if $\kappa = 0$, and that the order of convergence $2n+m$ is exact (unimprovable) here – that is, the quantities $\tilde{R}_n$ and $\tilde{G}_n$ have the same order of convergence with respect to $h$ as the quantities $\left| \omega_{n+1} - \overline{\omega} \right|^2 - \left| \omega_n - \overline{\omega} \right|^2$ and $\left| \omega_G - \overline{\omega} \right|^2 - \left| \omega_n - \overline{\omega} \right|^2$ (see Remark 5 above). The proof is based on establishing the following equivalences (as $h \to 0$) analogous to those given in Proposition 9:

1) for $\kappa > 0$ ($\lambda$ is arbitrary): $\tilde{R}_n \sim \dfrac{\kappa^2 p^{2n+8}}{6^{n+2}} h^{2n+4} (e^{-2\kappa p^2 t^2})$, $\tilde{G}_n \sim \dfrac{\kappa^2 p^{2n+8}}{6^{n+2}} h^{2n+4} (e^{-2\kappa p^2 t^2})$;

2) for $\kappa = 0$: $\tilde{R}_n \sim -\lambda^2 p^{2n+8} a_n h^{2n+6} t^2$, $\tilde{G}_n \sim -\lambda^2 p^{2n+8} a_n h^{2n+6} t^2$ (where $a_1 = 7/6480$, $a_n = 1/(15 \cdot 6^{n+2})$, $n \geq 2$).

Since these equivalences differ from those of Proposition 9 only by positive factors ($e^{-2\kappa p^2 t^2}$ or $t^2$) independent of $h$, the asymptotical behaviours of $\tilde{R}_n$ and $\tilde{G}_n$ are similar to those of $\left| \omega_{n+1} - \overline{\omega} \right|^2 - \left| \omega_n - \overline{\omega} \right|^2$ and $\left| \omega_G - \overline{\omega} \right|^2 - \left| \omega_n - \overline{\omega} \right|^2$ as $h \to 0$.

**Proof.** Let us show that $\displaystyle\lim_{h \to 0} \dfrac{\tilde{R}_n}{h^{2n+m}} = \dfrac{\kappa^2 p^{2n+8}}{6^{n+2}} e^{-2\kappa p^2 t^2}$ if $\kappa > 0$ and $\displaystyle\lim_{h \to 0} \dfrac{\tilde{R}_n}{h^{2n+m}} = -\lambda^2 p^{2n+8} a_n t^2$ if $\kappa = 0$. We have $\tilde{R}_n = \left| e^{\overline{\omega}t} \right|^2 \left( \left| e^{t(\omega_{n+1} - \overline{\omega})} - 1 \right|^2 - \left| e^{t(\omega_n - \overline{\omega})} - 1 \right|^2 \right) = e^{-2\kappa p^2 t} \left( \left| e^{b+a} - 1 \right|^2 - \left| e^b - 1 \right|^2 \right)$, where $b \equiv t(\omega_n - \overline{\omega}) \xrightarrow{h \to 0} 0$, $a \equiv t(\omega_{n+1} - \omega_n) = t(-A_{n+1} - iB_{n+1}) \xrightarrow{h \to 0} 0$ (here $A_n$ and $B_n$ are defined as in the proof of Lemma 1). Using the asymptotical relation $e^a = 1 + \tilde{A}$ (where we denoted $\tilde{A} \equiv a + a^2/2! + \cdots$, $\tilde{A} = a + O(a^2)$ as $a \to 0$) and $|l_1 + l_2|^2 = |l_1|^2 + |l_2|^2 + 2(\operatorname{Re} l_1 \operatorname{Re} l_2 + \operatorname{Im} l_1 \operatorname{Im} l_2)$, which is valid for arbitrary complex numbers $l_1$ and $l_2$, we obtain the following representation:

$$\left| e^{b+a} - 1 \right|^2 - \left| e^b - 1 \right|^2 = \left| e^b - 1 + e^b \tilde{A} \right|^2 - \left| e^b - 1 \right|^2 = \left| e^b \tilde{A} \right|^2 + 2 \operatorname{Re}(e^b - 1) \operatorname{Re}(e^b \tilde{A}) + 2 \operatorname{Im}(e^b - 1) \operatorname{Im}(e^b \tilde{A}). \quad (*)$$

Note that $2n+m$ can be rewritten as $(2+\tilde{m}) + (2n+2)$, where $\tilde{m} = 0$ if $\kappa > 0$ and $\tilde{m} = 2$ if $\kappa = 0$. Then a direct arithmetical calculation of limits shows that $\displaystyle\lim_{h \to 0} \dfrac{e^b \tilde{A}}{h^{2n+2}} = -\left( \dfrac{\kappa p^{2n+4}}{6^{n+1}} + i \dfrac{\lambda p^{2n+3}}{6^{n+1}} \right) t$,

$$\lim_{h \to 0} \dfrac{\left| e^b \tilde{A} \right|^2}{h^{2n+m}} = \begin{cases} \dfrac{\lambda^2 p^{10} t^2}{1296} & \text{if } \kappa = 0, \ n = 1, \\ 0 & \text{otherwise,} \end{cases} \qquad \lim_{h \to 0} \dfrac{e^b - 1}{h^{2+\tilde{m}}} = \begin{cases} \dfrac{-\kappa p^4 t}{12} & \text{if } \kappa > 0, \\ i \dfrac{\lambda p^5 t}{c_n} & \text{if } \kappa = 0, \end{cases}$$

where $c_1 = 30$ and $c_n = 180$ ($n \geq 2$) in the last expression. Combining these results and (*), we obtain

$$\lim_{h \to 0} \frac{\left|e^{b+a}-1\right|^2 - \left|e^b-1\right|^2}{h^{2n+m}} = \frac{\kappa^2 p^{2n+8}}{6^{n+2}} t^2 \text{ if } \kappa > 0 \text{ and } \lim_{h \to 0} \frac{\left|e^{b+a}-1\right|^2 - \left|e^b-1\right|^2}{h^{2n+m}} = -\lambda^2 p^{2n+8} a_n t^2 \text{ if } \kappa = 0,$$

which immediately implies the desired result about $\lim_{h \to 0} \frac{\tilde{R}_n}{h^{2n+m}}$.

Similarly, defining $a$ and $b$ as $b \equiv t(\omega_n - \overline{\omega}) \xrightarrow[h \to 0]{} 0$, $a \equiv t(\omega_G - \omega_n) \xrightarrow[h \to 0]{} 0$ and making analogous calculations, we get the corresponding results for the limit $\lim_{h \to 0} \frac{\tilde{G}_n}{h^{2n+m}}$:

$$\lim_{h \to 0} \frac{\tilde{G}_n}{h^{2n+m}} = \frac{\kappa^2 p^{2n+8}}{6^{n+2}} e^{-2\kappa p^2 t} t^2 \text{ if } \kappa > 0 \text{ and } \lim_{h \to 0} \frac{\tilde{G}_n}{h^{2n+m}} = -\lambda^2 p^{2n+8} a_n t^2 \text{ if } \kappa = 0. \blacksquare$$

**Remark 8.** In the text above we assumed everywhere that $p \neq 0$. Note that the case $p = 0$ is trivial, since $\omega_n = \omega_G = \omega_L = \overline{\omega} = 0$ in this case.

Let us summarize the results of the current Section. We showed that increasing the number of corrections leads to error increase in the presence of diffusion terms (see Proposition 3). We also showed that all the corrected schemes are more accurate than the consistent Galerkin formulation for problems with diffusion (see Proposition 5). For the pure (diffusionless) transport problems the situation is completely opposite – i.e., increasing the number of corrections should improve the accuracy of the numerical solution (see Proposition 4), and the consistent Galerkin formulation produces more accurate results than all the corrected schemes (see Proposition 6). We discussed that the first correction (as the most frequently used case of application of the Guermond-Pasquetti technique in literature) always provides better accuracy than the lumped Galerkin formulation (see Proposition 7). We also investigated the differences between the consistent solution and the corrected ones, and showed that increasing the number of corrections makes solutions of the corrected schemes closer to the consistent solution in all cases (see Proposition 8). We also derived asymptotical formulas (as $h \to 0$) for expressions $\left|\omega_{n+1} - \overline{\omega}\right|^2 - \left|\omega_n - \overline{\omega}\right|^2$, $\left|\omega_G - \overline{\omega}\right|^2 - \left|\omega_n - \overline{\omega}\right|^2$, $\left\|e^{t\omega_{n+1}} e^{ikhp} - e^{\overline{\omega}t} e^{ikhp}\right\|_\infty^2 - \left\|e^{t\omega_n} e^{ikhp} - e^{\overline{\omega}t} e^{ikhp}\right\|_\infty^2$ and $\left\|e^{t\omega_G} e^{ikhp} - e^{\overline{\omega}t} e^{ikhp}\right\|_\infty^2 - \left\|e^{t\omega_n} e^{ikhp} - e^{\overline{\omega}t} e^{ikhp}\right\|_\infty^2$ (see Remarks 5-7) that will be useful in verifying numerical calculations below and comparing theoretical orders of convergence with empirical ones.

### 3. Numerical examples

Let us compare the accuracy of the corrected and consistent schemes for the classical Galerkin FEM to confirm the theoretical results obtained. As in the paper [11], the time stepping is done with the standard explicit forth-order Runge-Kutta method (RK4) (with time step $\tau = 10^{-11}$, $\tau = 10^{-8}$, $\tau = 10^{-7}$ for examples 1-2, 3-4 and 5-7, respectively). This method and very small time steps were used to minimize the effect of discretization on the time variable, thus to ensure that the time error contributed by the time approximation is negligible in comparison with the error induced by spatial approximations (see [11]). It should be noted that no spurious oscillations appeared in all the examples considered below. The initial condition (for $t = 0$) and the boundary conditions are determined from known analytical solutions by their continuous extensions to the corresponding bounding hyperplanes (but for one-dimensional problems we use periodic boundary conditions). In all one-dimensional examples we use uniform meshes with the step $h$ to validate the theoretical estimates of Fourier analysis (Subsection 2.3). In particular, we determine the empirical orders $s$ of convergence in estimates of the form $O(h^s)$ for the differences of squared absolute errors for the corrected, consistent and analytical solutions (see Remarks 7 and 5), and reveal that these empirical orders converge to the corresponding theoretical ones established in Remark 7.

In all two-dimensional (three-dimensional) examples we use linear triangular (tetrahedral) 3-noded (4-noded) Lagrange-type elements obtained with a Delaunay triangulation algorithm.

Calculations in all the examples considered below show that in the presence of diffusion terms, when mesh size tends to zero, the first corrected scheme gives the most accurate results and increasing the order of correction only worsens the accuracy; finally, the least accurate results are given by the original Galerkin formulation with the consistent mass matrix in this case. These calculations confirm the theoretical conclusions of Fourier analysis from Subsection 2.3 established for one-dimensional harmonics.

*One-dimensional examples*

**Example 1.** Consider the initial-boundary value problem for the equation (1) with the known harmonic solution $u_{\text{exact}}(t,x) = Ae^{\overline{\omega}t}e^{ixp}$, where $\lambda = 1$, $\kappa = 10^{-2}$, $A = 1$, $p = 3\pi$, $x \in [0, 10]$.

Here and in what follows we denote the difference $u_{\text{numerical}}(t,x) - u_{\text{exact}}(t,x)$ by $err_n(t,x)$, $err_G(t,x)$ and $err_L(t,x)$ for the $n$-th corrected scheme, consistent Galerkin scheme (Eq. (2) for 1D case) and lumped scheme (3), respectively; then $\|err_j\|_\infty \equiv \max_{x_k}|err_j(t,x_k)|$ and $\|err_j\|_{\infty,\text{rel}} \equiv \max_{x_k}\left(\dfrac{|err_j(t,x_k)|}{|u_{\text{exact}}(t,x_k)|}\right)$ denote the maximum norm of the absolute and relative errors over all the grid points $x_k$ in the mesh, respectively. It should be noted that in Examples 4 and 7 below we omitted the nodes $x_k$ with $|u_{\text{exact}}(t,x_k)| < 10^{-10}$ in the definition of $\|err_j\|_{\infty,\text{rel}}$ (inside the maximum) to avoid overflow errors; other examples did not contain such points.

Figure 3 contains the plots of the functions $f_j(z)$ and $g_j(z)$ ($1 \leq j \leq 4$) for this example. The functions $f_1(z)$, $f_2(z)$, $g_1(z)$, $g_2(z)$ and $g_3(z)$ have the smallest positive roots at the points $c_1 \approx 0.194800$, $c_2 \approx 0.355856$, $c_3 \approx 0.195242$, $c_4 \approx 0.356028$ and $c_5 \approx 0.364753$, respectively, where these functions change their sign from plus to minus (note that we do not consider the function $f_3(z)$ that involves the 4-th correction, since calculations were made only for the first three corrected schemes). Due to the theoretical results of Subsection 2.3 (see Lemmas 1 and 2), the inequalities $|\omega_2 - \overline{\omega}| > |\omega_1 - \overline{\omega}|$, $|\omega_3 - \overline{\omega}| > |\omega_2 - \overline{\omega}|$, $|\omega_G - \overline{\omega}| > |\omega_1 - \overline{\omega}|$, $|\omega_G - \overline{\omega}| > |\omega_2 - \overline{\omega}|$ and $|\omega_G - \overline{\omega}| > |\omega_3 - \overline{\omega}|$ are valid if $h < c_1/|p|$, $h < c_2/|p|$, $h < c_3/|p|$, $h < c_4/|p|$ and $h < c_5/|p|$, respectively (note that $|p| = p = 3\pi$ in this example); the last inequalities are equivalent to the inequalities $N \geq 485$, $N \geq 266$, $N \geq 484$, $N \geq 266$ and $N \geq 260$, respectively, where $N$ is the total number of nodes (note that $h = (10-0)/(N-1)$ in this example). The corresponding calculations are presented in Table 1, where we also report the results for the numbers $N$ that are one less than the indicated "boundary" values, to illustrate the effect of violation of above inequalities.

*Table 1. Errors for Example 1, $t = 10^{-1}$.*

| Value | The number of spatial nodes, $N$ | | | | | | |
|---|---|---|---|---|---|---|---|
| | 259 | 260 | 265 | 266 | 483 | 484 | 485 |
| $\|err_1\|_{\infty,\text{rel}}$ | 1.0964e-03 | 1.0861e-03 | 1.0371e-03 | 1.0277e-03 | 2.8345e-04 | 2.8226e-04 | 2.8108e-04 |
| $\|err_2\|_{\infty,\text{rel}} - \|err_1\|_{\infty,\text{rel}}$ | -1.0017e-04 | -9.7696e-05 | -8.6255e-05 | -8.4145e-05 | -2.7375e-08 | -1.2436e-08 | 2.1708e-09 |
| $\|err_3\|_{\infty,\text{rel}} - \|err_2\|_{\infty,\text{rel}}$ | -5.3026e-08 | -4.4055e-08 | -5.8660e-09 | 5.7960e-10 | 1.6066e-08 | 1.5896e-08 | 1.5727e-08 |
| $\|err_G\|_{\infty,\text{rel}} - \|err_1\|_{\infty,\text{rel}}$ | -1.0023e-04 | -9.7740e-05 | -8.6260e-05 | -8.4143e-05 | -1.1206e-08 | 3.5619e-09 | 1.7999e-08 |
| $\|err_G\|_{\infty,\text{rel}} - \|err_2\|_{\infty,\text{rel}}$ | -5.3100e-08 | -4.3969e-08 | -5.1392e-09 | 1.4074e-09 | 1.6169e-08 | 1.5998e-08 | 1.5828e-08 |
| $\|err_G\|_{\infty,\text{rel}} - \|err_3\|_{\infty,\text{rel}}$ | -7.3585e-11 | 8.5912e-11 | 7.2678e-10 | 8.2777e-10 | 1.0336e-10 | 1.0183e-10 | 1.0033e-10 |
| $\|err_L\|_{\infty,\text{rel}}$ | 2.0855e-02 | 2.0695e-02 | 1.9924e-02 | 1.9774e-02 | 6.0018e-03 | 5.9770e-03 | 5.9524e-03 |

Here and in what follows we also denote by $P_{k,j}$ the value of $\ln(A_{h_{i-1}}(k,j)/A_{h_i}(k,j))/\ln(h_{i-1}/h_i)$ (here $h_i$ is the spatial step for the $i$-th column of the table, $A_{h_i}(k,j)$ is the value of $A(k,j) = \|err_k\|_\infty^2 - \|err_j\|_\infty^2$ for the step $h_i$) giving an empirical order of convergence with decreasing the spatial step, to confirm the results of Fourier analysis regarding the orders of convergence for differences between squared absolute errors of various numerical schemes (see Remark 7). Due to this definition $A(n+1,n)$ and $A(G,n)$ are equal to $\tilde{R}_n$ and $\tilde{G}_n$ from Remark 7 above, respectively; thus the empirical orders $P_{n+1,n}$ and $P_{G,n}$ (or $P_{n,n+1}$ and $P_{n,G}$ owing to symmetry in their definition) should converge to the number $2n+m$ (where $m=4$ if $\kappa > 0$ and $m=6$ if $\kappa = 0$) as $h \to 0$, as predicted by the theory developed in Subsection 2.3 (see Remarks 7 and 5).

The errors together with the corresponding empirical orders of convergence are given in Table 2 (note that in contrast to the calculations reported in Table 1 above, all the inequalities $|\omega_2 - \overline{\omega}| > |\omega_1 - \overline{\omega}|$, $|\omega_3 - \overline{\omega}| > |\omega_2 - \overline{\omega}|$, $|\omega_G - \overline{\omega}| > |\omega_1 - \overline{\omega}|$, $|\omega_G - \overline{\omega}| > |\omega_2 - \overline{\omega}|$ and $|\omega_G - \overline{\omega}| > |\omega_3 - \overline{\omega}|$ are now valid for all $N$ used in Table 2, since $N \geq 485$). It is clearly seen from Table 2 that the numbers $P_{2,1}$, $P_{3,2}$, $P_{G,1}$, $P_{G,2}$ and $P_{G,3}$ converge monotonically to 6, 8, 6, 8 and 10, respectively (that is, to the number $2n+4$ for $n=1$, $n=2$, $n=1$, $n=2$ and $n=3$, respectively).

*Table 2. Errors for Example 1, $t = 10^{-1}$.*

| Value | The number of spatial nodes, $N$ | | | | | | | | | |
|---|---|---|---|---|---|---|---|---|---|---|
| | 501 | 601 | 701 | 801 | 901 | 1001 | 1101 | 1201 | 1501 | 2501 |
| $\|err_1\|_{\infty,\text{rel}}$ | 2.6315e-04 | 1.8228e-04 | 1.3385e-04 | 1.0248e-04 | 8.0992e-05 | 6.5621e-05 | 5.4245e-05 | 4.5591e-05 | 2.9192e-05 | 1.0516e-05 |
| $\|err_2\|_{\infty,\text{rel}}$ | 2.6335e-04 | 1.8280e-04 | 1.3427e-04 | 1.0278e-04 | 8.1203e-05 | 6.5769e-05 | 5.4352e-05 | 4.5669e-05 | 2.9226e-05 | 1.0521e-05 |
| $\|err_3\|_{\infty,\text{rel}} -$ $- \|err_2\|_{\infty,\text{rel}}$ | 1.3287e-08 | 4.9787e-09 | 2.1007e-09 | 9.7956e-10 | 4.9561e-10 | 2.6811e-10 | 1.5331e-10 | 9.1848e-11 | 2.4519e-11 | 1.1674e-12 |
| $\|err_G\|_{\infty,\text{rel}} -$ $- \|err_3\|_{\infty,\text{rel}}$ | 7.9339e-11 | 2.0564e-11 | 6.3641e-12 | 2.2701e-12 | 9.0703e-13 | 3.9732e-13 | 1.8773e-13 | 9.4487e-14 | 1.6160e-14 | 2.7544e-16 |
| $\|err_L\|_{\infty,\text{rel}}$ | 5.5781e-03 | 3.8757e-03 | 2.8483e-03 | 2.1811e-03 | 1.7236e-03 | 1.3963e-03 | 1.1540e-03 | 9.6974e-04 | 6.2070e-04 | 2.2348e-04 |
| $P_{2,1}$ | — | -3.3706 | 3.3853 | 4.5353 | 5.0222 | 5.2878 | 5.4527 | 5.5637 | 5.6942 | 5.8581 |
| $P_{3,2}$ | — | 7.3865 | 7.5994 | 7.7147 | 7.7853 | 7.8321 | 7.8648 | 7.8886 | 7.9189 | 7.9605 |
| $P_{G,1}$ | — | -3.0604 | 3.4149 | 4.5483 | 5.0299 | 5.2930 | 5.4566 | 5.5667 | 5.6963 | 5.8591 |
| $P_{G,2}$ | — | 7.3966 | 7.6065 | 7.7200 | 7.7894 | 7.8354 | 7.8675 | 7.8909 | 7.9206 | 7.9613 |
| $P_{G,3}$ | — | 9.4082 | 9.6103 | 9.7212 | 9.7896 | 9.8351 | 9.8670 | 9.8904 | 9.9203 | 9.9687 |

Let us now consider the same problem but with $\kappa = 0$. For this case the errors together with the corresponding empirical orders of convergence are given in Table 3. The standard condition $|z| = |ph| \leq \pi$ is satisfied for the data reported in this table.

From Table 3, it is clearly seen that the reported numbers $P_{1,2}$, $P_{2,3}$, $P_{1,G}$, $P_{2,G}$ and $P_{3,G}$ converge monotonically to 8, 10, 8, 10 and 12, respectively (that is, to the number $2n+6$ for $n=1$, $n=2$, $n=1$, $n=2$ and $n=3$, respectively).

*Table 3. Errors for Example 1 (pure convection ), $t = 10^{-1}$.*

| Value | The number of spatial nodes, $N$ | | | | | | |
|---|---|---|---|---|---|---|---|
| | 501 | 601 | 701 | 801 | 901 | 1001 | 1101 |
| $\|err_L\|_{\infty,\text{rel}}$ | 5.5712e-03 | 3.8710e-03 | 2.8449e-03 | 2.1786e-03 | 1.7216e-03 | 1.3947e-03 | 1.1527e-03 |
| $\|err_1\|_{\infty,\text{rel}}$ | 3.9493e-05 | 1.9070e-05 | 1.0302e-05 | 6.0417e-06 | 3.7731e-06 | 2.4761e-06 | 1.6915e-06 |
| $\|err_2\|_{\infty,\text{rel}}$ | 6.8320e-06 | 3.2622e-06 | 1.7502e-06 | 1.0219e-06 | 6.3623e-07 | 4.1661e-07 | 2.8414e-07 |
| $\|err_3\|_{\infty,\text{rel}}$ | 6.6392e-06 | 3.1973e-06 | 1.7244e-06 | 1.0103e-06 | 6.3050e-07 | 4.1357e-07 | 2.8242e-07 |
| $\|err_3\|_{\infty,\text{rel}} -$ $- \|err_G\|_{\infty,\text{rel}}$ | 1.1453e-09 | 2.6733e-10 | 7.8059e-11 | 2.6860e-11 | 1.0479e-11 | 4.5139e-12 | 2.1068e-12 |
| $P_{1,2}$ | — | 7.9821 | 7.9872 | 7.9904 | 7.9926 | 7.9941 | 7.9952 |
| $P_{2,3}$ | — | 10.0065 | 10.0047 | 10.0036 | 10.0028 | 10.0023 | 10.0019 |
| $P_{1,G}$ | — | 7.9850 | 7.9893 | 7.9920 | 7.9938 | 7.9950 | 7.9959 |
| $P_{2,G}$ | — | 10.0161 | 10.0117 | 10.0088 | 10.0069 | 10.0055 | 10.0045 |
| $P_{3,G}$ | — | 11.9875 | 11.9910 | 11.9932 | 11.9947 | 11.9958 | 11.9966 |

**Example 2.** Let us consider the problem with the known harmonic solution $u_{\text{exact}}(t,x) = Ae^{\overline{\omega}t}e^{ixp}$, where $\lambda = 1$, $\kappa = 0$, $A = 1$, $p = 20\pi$ on the interval $x \in [0,1]$. For this case the errors together with the corresponding empirical orders of convergence are given in Table 4. The standard condition $|z| = |ph| \leq \pi$ is satisfied for the data reported in this table.

Again, from Table 4, it is clearly seen that the reported numbers $P_{1,2}$, $P_{2,3}$, $P_{1,G}$, $P_{2,G}$ and $P_{3,G}$

converge monotonically to 8, 10, 8, 10 and 12, respectively (that is, to the number $2n+6$ for $n=1$, $n=2$, $n=1$, $n=2$ and $n=3$, respectively).

*Table 4. Errors for Example 2, $t=1$.*

| Value | The number of spatial nodes, $N$ | | | | | | |
|---|---|---|---|---|---|---|---|
| | 501 | 601 | 701 | 801 | 901 | 1001 | 1101 |
| $\|err_L\|_{\infty,rel}$ | 1.6505e-01 | 1.1471e-01 | 8.4312e-02 | 6.4565e-02 | 5.1021e-02 | 4.1331e-02 | 3.4159e-02 |
| $\|err_1\|_{\infty,rel}$ | 5.2129e-04 | 2.5154e-04 | 1.3582e-04 | 7.9634e-05 | 4.9723e-05 | 3.2627e-05 | 2.2286e-05 |
| $\|err_2\|_{\infty,rel}$ | 8.8350e-05 | 4.2415e-05 | 2.2832e-05 | 1.3360e-05 | 8.3304e-06 | 5.4608e-06 | 3.7274e-06 |
| $\|err_3\|_{\infty,rel}$ | 8.7212e-05 | 4.2034e-05 | 2.2681e-05 | 1.3292e-05 | 8.2968e-06 | 5.4429e-06 | 3.7173e-06 |
| $\|err_3\|_{\infty,rel} - \|err_G\|_{\infty,rel}$ | 2.9994e-09 | 6.9850e-10 | 2.0370e-10 | 7.0034e-11 | 2.7305e-11 | 1.1756e-11 | 5.4850e-12 |
| $P_{1,2}$ | — | 7.9921 | 7.9943 | 7.9958 | 7.9967 | 7.9974 | 7.9978 |
| $P_{2,3}$ | — | 10.0030 | 10.0022 | 10.0016 | 10.0014 | 10.0007 | 10.0003 |
| $P_{1,G}$ | — | 7.9933 | 7.9952 | 7.9964 | 7.9972 | 7.9978 | 7.9982 |
| $P_{2,G}$ | — | 10.0074 | 10.0053 | 10.0040 | 10.0032 | 10.0021 | 10.0015 |
| $P_{3,G}$ | — | 11.9958 | 11.9964 | 11.9973 | 11.9983 | 11.9990 | 11.9997 |

The empirical orders $P_{i,j}$ (relative to $h$) of the differences of squared absolute errors indicated in Tables 2-4 fully correspond (converge with decreasing $h$) to the results of the Fourier analysis and the theoretical estimates obtained in Subsection 2.3 (see Remarks 7 and 5). Thus, the examples considered confirm the conclusions of Fourier-analysis regarding the accuracy of all the considered semi-discrete schemes with decreasing $h$.

*Two-dimensional and three-dimensional examples*

**Example 3.** Let us consider the problem for the two-dimensional convection-diffusion equation $u_t + \lambda_1 u_x + \lambda_2 u_y = \kappa(u_{xx} + u_{yy})$, $(x;y) \in [0,1] \times [0,1]$, with the known exact solution $u_{\text{exact}}(t,x,y) = A\exp(k_1 x + k_2 y + (k_1^2 + k_2^2)\kappa t)\exp\left(\frac{\lambda_1}{2\kappa}x + \frac{\lambda_2}{2\kappa}y - \frac{\lambda_1^2 + \lambda_2^2}{4\kappa}t\right)$, where $A=100$, $k_1=1$, $k_2=2$, $\lambda_1=1$, $\lambda_2=3/2$, $\kappa=1$. For this case the errors are given in Table 5. For further detailing we also reported the relative Euclidian error $\|err_j\|_{2,\text{dis}} \equiv \sqrt{\sum_{i=1}^{N}(err_j(t,\vec{x}_i))^2} / \sqrt{\sum_{i=1}^{N}(u_{\text{exact}}(t,\vec{x}_i))^2}$ (the sums are taken over all grid nodes $\{\vec{x}_i\}$). The mesh is assumed to be uniform along each direction.

*Table 5. Errors for Example 3, $t=1/2$.*

| Value | Counts of spatial nodes ($N_x$; $N_y$) | | | | | |
|---|---|---|---|---|---|---|
| | (15; 25) | (19; 29) | (25; 35) | (29; 39) | (35; 45) | (39; 49) |
| $\|err_1\|_{\infty,rel}$ | 6.5800e-4 | 4.2984e-4 | 2.6805e-4 | 2.0838e-4 | 1.5126e-4 | 1.2569e-4 |
| $\|err_2\|_{\infty,rel}$ | 6.6100e-4 | 4.3173e-4 | 2.6885e-4 | 2.0911e-4 | 1.5175e-4 | 1.2604e-4 |
| $\|err_3\|_{\infty,\text{rel}}$ | 6.6151e-4 | 4.3205e-4 | 2.6899e-4 | 2.0922e-4 | 1.5182e-4 | 1.2610e-4 |
| $\|err_4\|_{\infty,\text{rel}}$ | 6.6167e-4 | 4.3214e-4 | 2.6903e-4 | 2.0925e-4 | 1.5184e-4 | 1.2612e-4 |
| $\|err_G\|_{\infty,\text{rel}}$ | 6.6176e-4 | 4.3220e-4 | 2.6906e-4 | 2.0928e-4 | 1.5185e-4 | 1.2613e-4 |
| $\|err_1\|_{2,\text{dis}}$ | 1.8062e-4 | 1.2098e-4 | 7.7234e-5 | 6.0699e-5 | 4.4649e-5 | 3.7348e-5 |
| $\|err_2\|_{2,\text{dis}}$ | 1.8615e-4 | 1.2436e-4 | 7.9097e-5 | 6.2033e-5 | 4.5512e-5 | 3.8016e-5 |
| $\|err_3\|_{2,\text{dis}}$ | 1.8707e-4 | 1.2489e-4 | 7.9375e-5 | 6.2227e-5 | 4.5635e-5 | 3.8110e-5 |
| $\|err_4\|_{2,\text{dis}}$ | 1.8734e-4 | 1.2504e-4 | 7.9452e-5 | 6.2281e-5 | 4.5669e-5 | 3.8136e-5 |
| $\|err_G\|_{2,\text{dis}}$ | 1.8752e-4 | 1.2514e-4 | 7.9504e-5 | 6.2317e-5 | 4.5692e-5 | 3.8154e-5 |

**Example 4.** Let us consider the two-dimensional pure transport problem: $u_t + \lambda_1 u_x + \lambda_2 u_y = 0$, $(x;y) \in [0,1] \times [0,1]$, with the exact solution $u_{\text{exact}}(t,x,y) = \cos(2\pi(x-\lambda_1 t))\cos(2\pi(y-\lambda_2 t))$, where $\lambda_1 = 1$, $\lambda_2 = 3/2$. As in previous example, the mesh is assumed to be uniform along each direction. For this case the errors are given in Table 6.

*Table 6. Errors for Example 4, $t = 1/2$.*

| Value | Counts of spatial nodes ($N_x$; $N_y$) | | | | | |
|---|---|---|---|---|---|---|
| | (15; 25) | (19; 29) | (25; 35) | (29; 39) | (35; 45) | (39; 49) |
| $\|err_1\|_{\infty,\text{rel}}$ | 8.7340e-1 | 7.0159e-1 | 1.5334e-1 | 2.6956e-1 | 4.6210e-1 | 4.5959e-1 |
| $\|err_2\|_{\infty,\text{rel}}$ | 2.9288e-1 | 2.5134e-1 | 4.8985e-2 | 8.4199e-2 | 1.3308e-1 | 1.1962e-1 |
| $\|err_3\|_{\infty,\text{rel}}$ | 8.4041e-2 | 7.2884e-2 | 3.8187e-2 | 2.6194e-2 | 6.6196e-2 | 7.2150e-2 |
| $\|err_4\|_{\infty,\text{rel}}$ | 5.4847e-2 | 4.9143e-2 | 2.6638e-2 | 1.8352e-2 | 4.2513e-2 | 4.5909e-2 |
| $\|err_G\|_{\infty,\text{rel}}$ | 1.4612e-2 | 6.8652e-3 | 1.7686e-3 | 1.9321e-3 | 2.0885e-3 | 2.0157e-3 |
| $\|err_1\|_{2,\text{dis}}$ | 1.9701e-2 | 1.0722e-2 | 5.4253e-3 | 3.7973e-3 | 2.4469e-3 | 1.9118e-3 |
| $\|err_2\|_{2,\text{dis}}$ | 7.3973e-3 | 4.1823e-3 | 2.2105e-3 | 1.5821e-3 | 1.0449e-3 | 8.2651e-4 |
| $\|err_3\|_{2,\text{dis}}$ | 3.7435e-3 | 2.1309e-3 | 1.1370e-3 | 8.1739e-4 | 5.4249e-4 | 4.3003e-4 |
| $\|err_4\|_{2,\text{dis}}$ | 2.1812e-3 | 1.2109e-3 | 6.3944e-4 | 4.5817e-4 | 3.0294e-4 | 2.3966e-4 |
| $\|err_G\|_{2,\text{dis}}$ | 7.9943e-4 | 2.8098e-4 | 8.9348e-5 | 5.0070e-5 | 2.4931e-5 | 1.6911e-5 |

**Example 5.** Let us consider the problem for the three-dimensional convection-diffusion equation $u_t + \lambda_1 u_x + \lambda_2 u_y + \lambda_3 u_z = \kappa(u_{xx} + u_{yy} + u_{zz})$, $(x;y;z) \in [0,1]^3$, with the exact solution

$$u_{\text{exact}}(t,x,y,z) = A \exp\left(k_1 x + k_2 y + k_3 z + (k_1^2 + k_2^2 + k_3^2)\kappa t\right) \exp\left(\frac{\lambda_1}{2\kappa} x + \frac{\lambda_2}{2\kappa} y + \frac{\lambda_3}{2\kappa} z - \frac{\lambda_1^2 + \lambda_2^2 + \lambda_3^2}{4\kappa} t\right),$$

where $A = 10$, $k_1 = 1$, $k_2 = 3/2$, $k_3 = 2$, $\lambda_1 = 1$, $\lambda_2 = 3/2$, $\lambda_3 = 2$, $\kappa = 2$. At first we consider two different grids with the same number of nodes – a uniform $11 \times 13 \times 15$-noded grid (Fig. 5a) and an unstructured non-uniform grid where the internal nodes (i.e., nodes inside the domain) are randomly distributed (Fig. 5b). Subsequently, we applied the Delaunay triangulation algorithm to obtain finite-element decompositions in both cases. Finally, let us also consider this problem on $13\times15\times17$-noded grids – as in the previous case, on a uniform grid (Fig. 6a) and on a grid with randomly distributed internal nodes (Fig. 6b). Corresponding errors for the time $t=1/2$ are given in Table 7.

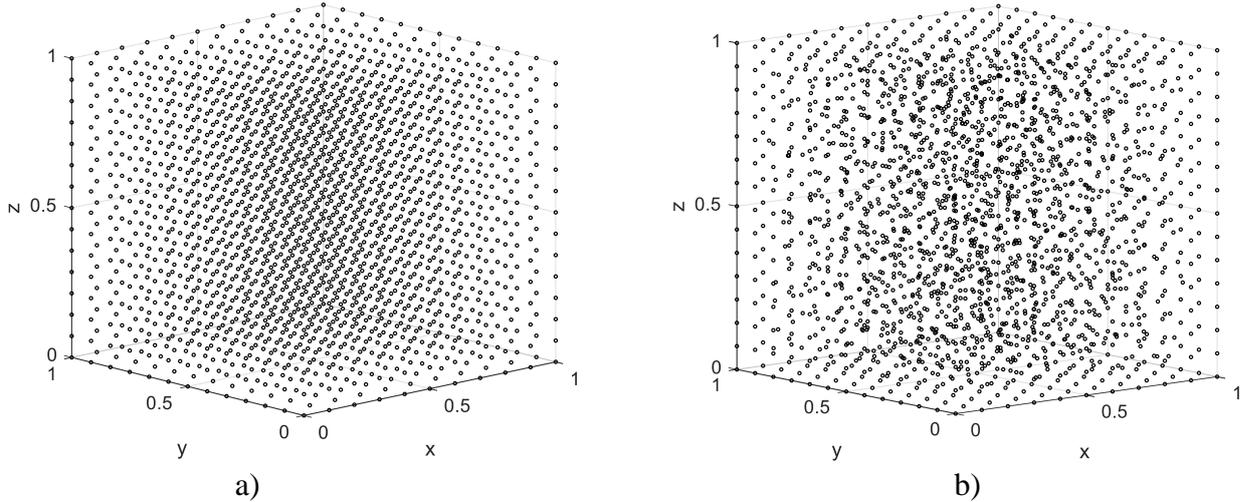

*Figure 5. 11x13x15-noded grids: a) uniform grid; b) grid with randomly distributed nodes.*

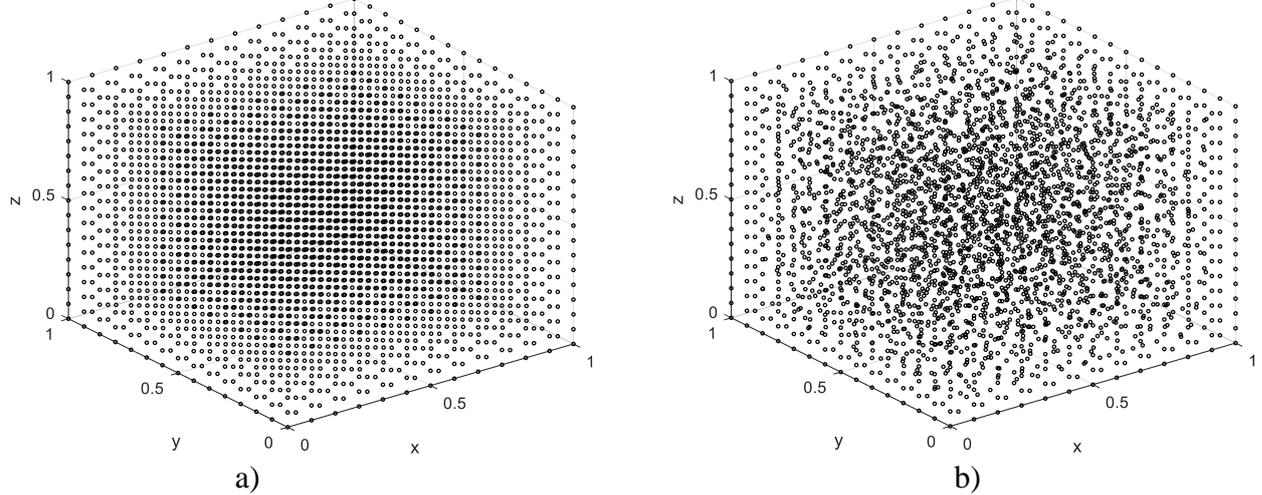

*Figure 6. 13x15x17-noded grids: a) uniform grid; b) grid with randomly distributed nodes.*

*Table 7. Errors for Example 5, $t=1/2$.*

| Value | $11\times13\times15$-noded uniform grid | $11\times13\times15$-noded random grid | $13\times15\times17$-noded uniform grid | $13\times15\times17$-noded random grid |
|---|---|---|---|---|
| $\|err_1\|_{\infty,\text{rel}}$ | 1.6885e-3 | 2.2759e-2 | 1.3321e-3 | 1.6624e-2 |
| $\|err_2\|_{\infty,\text{rel}}$ | 1.7062e-3 | 2.2863e-2 | 1.3449e-3 | 1.6626e-2 |
| $\|err_3\|_{\infty,\text{rel}}$ | 1.7087e-3 | 2.2900e-2 | 1.3465e-3 | 1.6629e-2 |
| $\|err_4\|_{\infty,\text{rel}}$ | 1.7090e-3 | 2.2919e-2 | 1.3469e-3 | 1.6632e-2 |
| $\|err_G\|_{\infty,\text{rel}}$ | 1.7094e-3 | 2.2943e-2 | 1.3473e-3 | 1.6637e-2 |
| $\|err_1\|_{2,\text{dis}}$ | 3.0321e-4 | 4.4343e-3 | 2.4788e-4 | 3.0021e-3 |
| $\|err_2\|_{2,\text{dis}}$ | 3.1808e-4 | 4.4602e-3 | 2.5951e-4 | 3.0209e-3 |
| $\|err_3\|_{2,\text{dis}}$ | 3.2098e-4 | 4.4666e-3 | 2.6164e-4 | 3.0252e-3 |
| $\|err_4\|_{2,\text{dis}}$ | 3.2183e-4 | 4.4688e-3 | 2.6224e-4 | 3.0268e-3 |
| $\|err_G\|_{2,\text{dis}}$ | 3.2239e-4 | 4.4708e-3 | 2.6265e-4 | 3.0285e-3 |

**Example 6.** Let us consider the problem for the three-dimensional convection-diffusion equation

$$\frac{\partial u}{\partial t} + \frac{\lambda}{t+1}\left(x\frac{\partial u}{\partial x} + y\frac{\partial u}{\partial y} + z\frac{\partial u}{\partial z}\right) = \kappa\left(\frac{\partial^2 u}{\partial x^2} + \frac{\partial^2 u}{\partial y^2} + \frac{\partial^2 u}{\partial z^2}\right), \quad (x; y; z) \in [0, 1]^3, \quad \text{with the exact solution}$$

$u_{\text{exact}}(t, x, y, z) = A \cdot (1 - 2\lambda)(x^2 + y^2 + z^2)(t+1)^{-2\lambda} + 6A\kappa(t+1)^{-2\lambda+1}$, where $A = 10^2$, $\kappa = 0.35$, $\lambda = 1$. We consider this problem on $11 \times 13 \times 15$-noded and $13 \times 15 \times 17$-noded grids from the previous Example 5 (see Figures 5-6). Corresponding errors for this case (for $t = 1/2$) are reported in Table 8.

*Table 8. Errors for Example 6, $t = 1/2$.*

| Value | $11 \times 13 \times 15$-noded uniform grid | $11 \times 13 \times 15$-noded random grid | $13 \times 15 \times 17$-noded uniform grid | $13 \times 15 \times 17$-noded random grid |
|---|---|---|---|---|
| $\|err_1\|_{\infty,\text{rel}}$ | 4.3634e-4 | 1.4331e-3 | 3.2482e-4 | 1.3372e-3 |
| $\|err_2\|_{\infty,\text{rel}}$ | 4.8768e-4 | 1.4511e-3 | 3.5697e-4 | 1.3645e-3 |
| $\|err_3\|_{\infty,\text{rel}}$ | 4.9362e-4 | 1.4564e-3 | 3.6059e-4 | 1.3702e-3 |
| $\|err_4\|_{\infty,\text{rel}}$ | 4.9466e-4 | 1.4580e-3 | 3.6137e-4 | 1.3722e-3 |
| $\|err_G\|_{\infty,\text{rel}}$ | 4.9566e-4 | 1.4594e-3 | 3.6203e-4 | 1.3746e-3 |
| $\|err_1\|_{2,\text{dis}}$ | 1.6658e-4 | 2.3380e-4 | 1.2473e-4 | 1.6989e-4 |
| $\|err_2\|_{2,\text{dis}}$ | 1.8461e-4 | 2.4153e-4 | 1.3697e-4 | 1.7446e-4 |
| $\|err_3\|_{2,\text{dis}}$ | 1.8764e-4 | 2.4344e-4 | 1.3887e-4 | 1.7564e-4 |
| $\|err_4\|_{2,\text{dis}}$ | 1.8845e-4 | 2.4422e-4 | 1.3938e-4 | 1.7615e-4 |
| $\|err_G\|_{2,\text{dis}}$ | 1.8898e-4 | 2.4525e-4 | 1.3971e-4 | 1.7684e-4 |

**Example 7.** Let us consider the three-dimensional pure transport problem:
$u_t + \lambda_1 u_x + \lambda_2 u_y + \lambda_3 u_z = 0$, $(x; y; z) \in [0, 1]^3$, with the exact solution $u_{\text{exact}}(t, x, y, z) = \sin(2\pi(x - \lambda_1 t))\sin(2\pi(y - \lambda_2 t))\sin(2\pi(z - \lambda_3 t))$, where $\lambda_1 = 1$, $\lambda_2 = -2$, $\lambda_3 = 3$. We consider this problem on $11 \times 13 \times 15$-noded and $13 \times 15 \times 17$-noded grids from the Example 5 (see Figures 5-6). Corresponding errors for this case are reported in Table 9.

*Table 9. Errors for Example 7, $t = 10^{-1}$.*

| Value | $11 \times 13 \times 15$-noded uniform grid | $11 \times 13 \times 15$-noded random grid | $13 \times 15 \times 17$-noded uniform grid | $13 \times 15 \times 17$-noded random grid |
|---|---|---|---|---|
| $\|err_1\|_{\infty,\text{rel}}$ | 4.7517e+0 | 9.7705e+3 | 3.1666e+1 | 7.1027e+4 |
| $\|err_2\|_{\infty,\text{rel}}$ | 1.8211e+0 | 6.9139e+3 | 9.6154e+0 | 5.9332e+4 |
| $\|err_3\|_{\infty,\text{rel}}$ | 7.5371e-1 | 5.4892e+3 | 3.8590e+0 | 4.8175e+4 |
| $\|err_G\|_{\infty,\text{rel}}$ | 3.7610e-1 | 3.2278e+3 | 3.0193e+0 | 1.2770e+4 |
| $\|err_1\|_{2,\text{dis}}$ | 6.2910e-2 | 1.5593e-1 | 4.4165e-2 | 1.3695e-1 |
| $\|err_2\|_{2,\text{dis}}$ | 2.6841e-2 | 1.4608e-1 | 1.8215e-2 | 1.2747e-1 |
| $\|err_3\|_{2,\text{dis}}$ | 1.5341e-2 | 1.4393e-1 | 9.9637e-3 | 1.2440e-1 |
| $\|err_G\|_{2,\text{dis}}$ | 8.0797e-3 | 1.4328e-1 | 4.2542e-3 | 1.2198e-1 |

As in all previous examples of pure transport problems, increasing the number of corrections improves the accuracy of the numerical solution and the consistent Galerkin formulation produces the most accurate results.

## 4. Conclusions

The paper provides a Fourier analysis of the Guermond-Pasquetti technique in application to pure transport and convection-diffusion problems. We show that increasing the number of corrections leads to error increase in the presence of diffusion terms (see Proposition 3). We also show that all the corrected schemes are more accurate than the consistent Galerkin formulation for problems with diffusion (see Proposition 5). For the pure (diffusionless) transport problems the situation is completely opposite – i.e., increasing the number of corrections should improve the accuracy of the numerical solution (see Proposition 4), and the consistent Galerkin formulation produces more accurate results than all the corrected schemes (see Proposition 6). We also investigate the differences between the consistent solution and the corrected ones, and show that increasing the number of corrections makes solutions of the corrected schemes closer to the consistent solution in all cases (see Proposition 8). Numerical examples confirmed the corresponding theoretical consequences of Fourier analysis.

**Acknowledgements.**
The author was partially supported by H2020-MSCA-RISE-2014 Project AMMODIT (Project number 645672). The author would like to thank the anonymous referees for their careful reading and helpful remarks. Also, the author expresses his deep gratitude to Prof. Alberto Redaelli and Dr. Filippo Piatti (Politecnico di Milano, Italy) for fruitful discussions and their warm welcome during several author's secondments. Above all, the author would like to thank his beloved wife Anna and his mother Liubov for their continual encouragement.